\documentclass{article}
\usepackage{amsfonts}
\usepackage{mathrsfs}

\usepackage{amsmath}
\usepackage{latexsym,amsfonts,amssymb}

\newcommand{\bfn}{{\bf{n}}}

\newcommand{\nec}{{\mbox{\rm\footnotesize nec}}}
\newcommand{\liv}{{\mbox{\rm\footnotesize liv}}}

\begin{document}

\title{Analysis of a free boundary problem modeling the growth of necrotic tumors\thanks{This work is supported by
  the National Natural Science Foundation of China under grant numbers 11571381.}}
\author{Shangbin Cui\\[0.2cm]
  {\small School of Mathematics, Sun Yat-Sen University, Guangzhou, Guangdong 510275,}\\
  {\small People's Republic of China. E-mail:\,cuishb@mail.sysu.edu.cn}}
\date{}
 \maketitle

\begin{abstract}
  In this paper we make rigorous analysis to a free boundary problem modeling the growth of a necrotic tumor.
  A remarkable feature of this free boundary problem is that it contains two different-type free surfaces: One is the
  tumor surface whose evolution is governed by an evolution equation and the other is the interface between the living
  shell of the tumor and the tumor's necrotic core which is an obstacle-type free surface, i.e., its evolution is not
  governed by an evolution equation but instead is determined by some stationary-type equation. In mathematics, the
  inner free surface is induced by discontinuity of the nonlinear reaction functions in this model, which causes the
  main difficulty of the analysis. Previous work on this model studies spherically symmetric situation which is in
  essence an one-dimension free boundary problem. The purpose of this paper is to make rigorous analysis in general
  spherically asymmetric situation. By applying the Nash-Moser implicit function theorem, we prove that the inner free
  surface is smooth and depends on the outer free surface smoothly when the outer free surface is a small perturbation
  of a sphere. By applying this result and some abstract theory of parabolic differential equations in Banach manifolds
  we prove that the unique radial stationary solution of this free boundary problem is asymptotically stable under
  small non-radial perturbation. 
\medskip

\medskip

\end{abstract}

\section{Introduction}
\setcounter{equation}{0}

\hskip 2em
  The study of mathematical theory of tumor growth has developed for over eighty years. It is motivated by a basic
  observation that under constant conditions, an evolutionary tumor undergoes three typical periods of growth: Firstly
  it goes through a nearly exponential growth period, next it experiences a nearly linear growth period, and finally
  it evolves into a stationary or dormant state \cite{AdmBel, Casey,Laird1,Laird2,Simp}. In the dormant state, the
  tumor usually contains an inner necrotic core made by dead cells, an outer proliferating shell occupied by proliferating
  cells, and an intermediate region occupied by quiescent living cells \cite{AdmBel,ByrC2,Green1, Green2}. During 1970's,
  Greenspan proposed the first mathematical model in the form of free boundary problem of reaction diffusion equations
  to explain this phenomenon \cite{Green1, Green2}. His model was very well improved by Byrne and Chaplain during 1990's
  \cite{ByrC1, ByrC2}. Since then many different tumor models have been established by different groups of researchers,
  cf. the reviewing articles \cite{tumrev1,Low} and references cited therein. Rigorous mathematical analysis of such
  free boundary problems has attracted much attention during the past twenty years and many interesting results have been
  obtained, cf. \cite{Cui1,Cui2, Cui3, CuiEsc1,CuiEsc2, CuiFri, FH1, FH2, FriRei1, FriRei2, WuCui, ZhuCui} and references
  cited therein.

  In this paper we study the following free boundary problem modeling the growth of a necrotic tumor:
\begin{equation}
\left\{
\begin{array}{rll}
   \Delta\sigma=&f(\sigma) &\quad\;\; \mbox{in} \;\; \Omega(t),\;\; t>0,\\
   -\Delta p=&g(\sigma)   &\quad\;\; \mbox{in} \;\; \Omega(t),\;\; t>0,\\
   \sigma=&\bar{\sigma}   &\quad\;\; \mbox{on} \;\; \partial\Omega(t),\;\; t>0,\\
   p=&\gamma\kappa  &\quad\;\; \mbox{on} \;\; \partial\Omega(t),\;\; t>0,\\
   V_n=&-\partial_{n}p  &\quad\;\; \mbox{on} \;\; \partial\Omega(t),\;\; t>0,\\
   \Omega(0)=&\Omega_0. &
\end{array}
\right.
\end{equation}
  Here $\Omega(t)$ is the domain in $\mathbb{R}^3$ occupied by the tumor at time $t$, $\sigma=\sigma(x,t)$ and $p=p(x,t)$
  are the nutrient concentration in the tumor region and the pressure between tumor cells, respectively, $\bar{\sigma}$
  is a positive constant reflecting the constant nutrient supply that the tumor receives from its host tissue, $\gamma$
  is a positive constant reflecting the surface tension of the tumor surface and is usually referred as {\em surface
  tension coefficient}, $\kappa$ is the mean curvature of the tumor surface $\partial\Omega(t)$ whose sign is designated
  by the convention that for the sphere it is positive, $V_{n}$ is the normal velocity of the tumor surface movement,
  $\partial_{n}$ represents the derivative in the direction of the outward normal $n$ of the tumor surface $\partial\Omega(t)$,
  $\Omega_{0}$ is the domain that the tumor initially occupies, and $f$, $g$ are given functions respectively having the
  following forms:
\begin{equation}
   f(\sigma)=\lambda\sigma H(\sigma-\hat{\sigma}),\qquad
   g(\sigma)=\mu(\sigma-\tilde{\sigma})H(\sigma-\hat{\sigma})-\nu,
\end{equation}
  where $H$ is the Heaviside function: $H(s)=1$ for $s>0$ and $H(s)=0$ for $s\leqslant 0$, and $\lambda$, $\mu$, $\nu$,
  $\hat{\sigma}$ and $\tilde{\sigma}$ are positive constants, with $\lambda$ being the consumption rate coefficient of
  nutrient by tumor cells, $\mu$ the proliferation rate coefficient of tumor cells (=the birth rate of tumor cells that
  a unit amount of nutrient can sustain), $\nu$ the dissolution rate of dead cells, $\hat{\sigma}$ a threshold value of
  nutrient concentration to sustain tumor cells alive and proliferating, i.e., only in the region where $\sigma>\hat{\sigma}$
  tumor cells are alive and proliferating, and $\tilde{\sigma}=\hat{\sigma}-(\nu/\mu)$. We assume that $0<\hat{\sigma}<
  \bar{\sigma}$, $\nu<\mu\hat{\sigma}$ (so that $0<\tilde{\sigma}<\hat{\sigma}$) and, for simplicity of notations
$$
  \lambda=\bar{\sigma}=1,
$$
  which can always be achieved through rescaling. For more information concerning biological implication of the above model,
  we refer the reader to see the reference \cite{ByrC2}.

  If instead of (1.2) the functions $f$, $g$ are smooth monotone increasing functions in $[0,\infty)$ and satisfy the properties
  $f(0)=0$, $g(0)<0$ and $g(\infty)>0$, the problem (1.1) models the growth of a nonnecrotic tumor, cf. \cite{ByrC1}, which
  has been intensively studied during the past twenty years, cf. \cite{Cui1}, \cite{CuiEsc1}, \cite{CuiEsc2}, \cite{FH1},
  \cite{FH2}, \cite{FriRei1} and references therein. It was proved that there exists a threshold value $\gamma^*>0$ for the
  surface tension coefficient $\gamma$, such that if $\gamma>\gamma^*$ then the unique radial stationary solution is
  asymptotically stable module translations, whereas if $\gamma<\gamma^*$ then it is unstable. In the necrotic case, analysis
  of the problem (1.1) is much harder: In addition to the outer free boundary $\partial\Omega(t)$ whose evolution is governed
  by the equation $V_n=-\partial_{\nu}p$, discontinuity of the functions $f$, $g$ at $\sigma=\hat{\sigma}$ produces an inner
  free surface $\Gamma(t)$ dividing the domain $\Omega(t)$ into two disjoint regions, with the outer region
\begin{equation}
   \Omega_{\rm liv}(t)=\{x\in\Omega(t):\sigma(x,t)>\hat{\sigma}\}
\end{equation}
  being the living shell and the inner region
\begin{equation}
   \Omega_{\rm nec}(t)={\rm int}\{x\in\Omega(t):\sigma(x,t)=\hat{\sigma}\}
\end{equation}
  the necrotic core of the tumor. $\Gamma(t)$ is therefore the common boundary of these two regions:
\begin{equation}
   \Gamma(t)=\partial\Omega_{\rm nec}(t)\cap\partial\Omega_{\rm liv}(t).
\end{equation}
  Main difficulty of analysis of the model (1.1) is caused by $\Gamma(t)$, which is implicitly contained in the problem (1.1)
  and there is not an obvious evolution equation to govern its movement.

  In \cite{Cui2} the spherically symmetric version of the problem (1.1) was studied, improving some earlier results of \cite{CuiFri}.
  It was proved that this problem has a unique radial stationary solution which is asymptotically stable under spherically symmetric
  perturbation. However, whether this stationary solution is asymptotically stable under spherically asymmetric perturbation
  has been kept unknown for over ten years. Not long ago some numerical results on this problem were obtained by Hao et al \cite{Hao1}.
  In this paper we aim at making a rigorous analysis to the problem (1.1) and establishing a similar result as that obtained in
  \cite{Cui1}, \cite{CuiEsc2}, \cite{FH1} for the nonnecrotic case, i.e.,
  we want to prove that there exists a threshold value $\gamma^*$ for the surface tension coefficient $\gamma$ such that if
  $\gamma>\gamma^*$ then the unique radial stationary solution ensured by the reference \cite{Cui2} is asymptotically stable module
  translations under spherically asymmetric perturbation, whereas if $\gamma<\gamma^*$ then it is unstable under spherically
  asymmetric perturbation. To get this result, the crucial step is to prove that the inner free interface $\Gamma(t)$ is smooth
  and depends on the outer free boundary $\partial\Omega(t)$ smoothly at least when the outer free boundary $\partial\Omega(t)$ is
  sufficiently close to the surface of a sphere. This will be proved by applying the Nash-Moser implicit function theorem. Note
  that this latter result has clearly its own independent significance.

  To state the main results of this paper, let us first recall some results obtained in \cite{Cui2}. Consider the following elliptic
  boundary value problem:
\begin{equation}
\left\{
\begin{array}{rll}
   \Delta\sigma=&f(\sigma) &\quad\;\; \mbox{in} \;\; \Omega,\\
   \sigma=&1   &\quad\;\; \mbox{on} \;\; \partial\Omega,
\end{array}
\right.
\end{equation}
  where $\Omega$ is a given bounded domain in ${\mathbf{R}}^3$ with a $C^2$ boundary and $f$ is as in (1.2) (with $\lambda=1$).
  It is not hard to see that as far as strong solution is concerned, the above problem is equivalent to the following obstacle
  problem:
$$
\left\{
\begin{array}{cl}
   -\Delta\sigma+\sigma\geqslant 0, \quad \sigma\geqslant\hat{\sigma}, \quad (-\Delta\sigma+\sigma)(\sigma-\hat{\sigma})=0 &\quad\;\,
   \mbox{in}\;\;\Omega,\\
   \sigma=1   &\quad \mbox{on} \;\; \partial\Omega.
\end{array}
\right.
$$
  It follows from the standard theory in obstacle problem (cf., e.g., \cite{Fried} and \cite{Tro}) that the problem (1.6) has
  a unique solution $\sigma\in\displaystyle\cap_{1\leqslant q<\infty}W^{2,q}(\Omega)$ which satisfies $\hat{\sigma}\leqslant
  \sigma\leqslant 1$ (see also Lemma 2.4 in Section 2). Besides, it is not very difficult to prove that (see Lemmas 3.1, 3.2
  of \cite{Cui2}) if $\Omega=B(0,R)$ for some $R>0$, then the unique solution of the above problem is a radial function, given
  by $\sigma(x)=U(r,R)$, where $r=|x|$, and $U(r,R)$ is defined as follows: Let $R^*$ be the unique positive number solving
  the equation
\begin{equation}
  \frac{\sinh\! R^*}{R^*}=\frac{1}{\hat{\sigma}}.
\end{equation}
  Then $U(r,R)=\displaystyle\frac{R\sinh r}{r\sinh R}$ for $0<R\leqslant R^*$ and
\begin{equation}
  U(r,R)=\left\{
\begin{array}{ll}
    \displaystyle\hat{\sigma}[\sinh(r- K)+ K\cosh(r- K)]/r &\quad\;\;
    \mbox{for} \;\;  K\leqslant r\leqslant R\\ [0.2cm]
   \;\;\hat{\sigma}   &\quad\;\; \mbox{for} \;\; r<K
\end{array}
\right.
\end{equation}
  for $R>R^*$, where for $R>R^*$, $K=K(R)$ is the unique solution of the following equation in the interval $(0,R)$:
\begin{equation}
  \sinh(R-K)+K\cosh(R-K)=\frac{R}{\hat{\sigma}}.
\end{equation}
  Now let
\begin{equation}
\begin{array}{rl}
  F(R)=&\displaystyle\frac{1}{4\pi R^3}\int_{B(0,R)}g(U(|x|,R))dx\\ [0.2cm]
   =&\displaystyle\frac{1}{R^3}\left\{\mu\int_{U(r,R)>\hat{\sigma}}\!\!\!\big(U(r,R)\!-\!
   \hat{\sigma}\big)r^2dr\!-\!\frac{1}{3}\nu R^3\right\}\\ [0.3cm]
   =&\displaystyle\left\{
\begin{array}{l}
  \displaystyle\mu\left\{\frac{R\coth R-1}{R^2}\!-\!\frac{1}{3}\hat{\sigma}\right\}, \quad 0<R\leqslant R^*,
\\ [3mm]
  \displaystyle\mu\hat{\sigma}G(R)\!-\!\frac{1}{3}\nu\Big(\frac{K}{R}\Big)^3\!-\!\frac{1}{3}\mu\hat{\sigma}, \quad R>R^*,
\end{array}
\right.
\end{array}
\end{equation}
  where $K=K(R)$ is as before, and
$$
  G(R)=\frac{1}{R^3}\left\{(R-K)\cosh(R-K)+(RK-1)\sinh(R-K)\!+\!\frac{1}{3}K^3\right\}, \quad R>R^*.
$$
  In case $\Omega(t)=B(0,R(t))$, the problem (1.1) reduces into the following initial value problem for a first-order differential
  equation:
$$
\left\{
\begin{array}{l}
  R'(t)=R(t)F(R(t)), \quad t>0,
\\
  R(0)=R_0,
\end{array}
\right.
$$
  where $R_0>0$ is the number such that $\Omega_0=B(0,R_0)$. It was proved (see Lemma 4.2 of \cite{Cui2}) that $F$ is continuously
  differentiable in $(0,\infty)$, $F'(R)<0$ for  all $R>0$, and
$$
  \lim_{R\to 0^+}F(R)=\frac{\mu}{3}(1-\hat{\sigma})>0, \quad
  \lim_{R\to \infty}F(R)=-\frac{\nu}{3}<0.
$$
  Hence the function $F$ has a unique positive root $R_s$, and $F(R)>0$ for $0<R<R_s$, $F(R)<0$ for $R>R_s$. It follows that for
  any $R_0>0$, the solution of the above problem  exists for all $t\geqslant0$, and
$$
  \lim_{t\to\infty}R(t)=R_s.
$$
  Moreover, since the function $R\mapsto\displaystyle\frac{\sinh R}{R}\cdot\frac{R\coth R-1}{(R)^2}$ is strictly monotone
  increasing, converges to $\displaystyle\frac{1}{3}$ as $R\to 0^+$ and tends to $+\infty$ as $R\to+\infty$, it follows that
  the following relation is true:
$$
  \frac{R^*\coth R^*-1}{(R^*)^2}>\frac{1}{3}\hat{\sigma}.
$$
  From this relation it follows that $R_s>R^*$, which implies that the dormant or stationary tumor must have a necrotic core
  with radius $K_s=K(R_s)$. It was also proved in \cite{Cui2} that the necrotic core is formed at finite time.

  Throughout this paper we use the notation $(r,\omega)$, $r\geqslant 0$, $\omega\in\mathbb{S}^2$, to denote the polar
  coordinate of point $x$ in $\mathbb{R}^3$, i.e., $r=|x|$ and $\omega=x/|x|$ for $x\in\mathbb{R}^3\backslash\{0\}$, and
  $|0|=0$. Given $R>R^*$ and $\rho,\eta\in C^{2+\mu}({\mathbf{S}}^2)$ ($0<\mu<1$) with $\|\rho\|_{C^{2+\mu}({\mathbf{S}}^2)}$
  and $\|\eta\|_{C^{2+\mu}({\mathbf{S}}^2)}$ sufficiently small, we denote
$$
  \Omega_{\rho}=\{x\in{\mathbf{R}}^3: r<R[1+\rho(\omega)]\}, \quad
  D_{\rho,\eta}=\{x\in{\mathbf{R}}^3:  K[1+\eta(\omega)]<r<R[1+\rho(\omega)]\},
$$
$$
   S_{\rho}=\{x\in{\mathbf{R}}^3: r=R[1+\rho(\omega)]\}, \quad \mbox{and} \quad
   \Gamma_{\eta}=\{x\in{\mathbf{R}}^3: r= K[1+\eta(\omega)]\}.
$$
  It is not hard to prove (see the remark following Lemma 2.4 in Section 2) that given $\rho\in C^{2+\mu}({\mathbf{S}}^2)$
  ($0<\mu<1$) with $\|\rho\|_{C^{2+\mu}({\mathbf{S}}^2)}$ sufficiently small and letting $\Omega=\Omega_{\rho}$, solving
  the problem (1.6) is equivalent to seeking $(\sigma,\eta)$ such that it solves the following problem:
\begin{equation}
\left\{
\begin{array}{rll}
   \Delta\sigma=&\sigma &  \quad  \mbox{in}\;\; D_{\rho,\eta},\\
   \sigma=&1 &  \quad  \mbox{on}\;\; S_{\rho},\\
   \sigma=&\hat{\sigma} &  \quad  \mbox{on}\;\; \Gamma_{\eta},\\
    \partial_{r}\sigma=&0 & \quad  \mbox{on}\;\; \Gamma_{\eta},
\end{array}
\right.
\end{equation}
  where $\partial_{r}$ denotes the derivative in radial direction. We note that if $\|\eta\|_{C^1({\mathbf{S}}^2)}$
  is sufficiently small then at any point in $\Gamma_{\eta}$, the radial direction is not tangent to $\Gamma_{\eta}$,
  so that the boundary value condition $(1.11)_4$ is regular. Later on we shall also use the following abbreviations:
$$
  D=D_{0,0}=B(0,R)\backslash\overline{B(0,K)}, \quad S_0=\partial B(0,R), \quad \mbox{and} \quad  \Gamma_0=\partial B(0,K),
$$
  where $K=K(R)$.
\medskip

  As we mentioned earlier, the problem (1.6) has a unique solution $\sigma\in\displaystyle\cap_{1\leqslant q<\infty}W^{2,q}(\Omega)$.
  Hence if $R>R^*$ then the problem (1.11) has a unique solution. What we are concerned with is regularity of the free boundary
  $\Gamma_{\eta}$, or equivalently the regularity of the function $\eta$, and, more importantly, the regularity
  of the mapping $\rho\mapsto\eta$. This leads to the following result:
\medskip

  {\bf Theorem 1.1}\ \ {\em Let $R>R^*$ be given. Let the integer $m\geqslant 2$ and the number $0<\mu<1$ be fixed. Then there exists
  a constant $\delta>0$ such that for any $\rho\in C^{m+\mu}({\mathbf{S}}^2)$ with $\|\rho\|_{C^{m+\mu}({\mathbf{S}}^2)}<\delta$,
  the problem $(1.11)$ has a unique solution $(\sigma,\eta)$ with $\eta\in C^{\infty}({\mathbf{S}}^2)$ and $\sigma\in
  C^{m+\mu}(\overline{D}_{\rho,\eta})\cap C^{\infty}(\overline{D}_{\rho,\eta}\backslash S_{\rho})$, and the mapping $\rho\mapsto\eta$
  from the open set $\|\rho\|_{C^{m+\mu}({\mathbf{S}}^2)}<\delta$ in $C^{m+\mu}({\mathbf{S}}^2)$ to the Frech\'{e}t space
  $C^{\infty}({\mathbf{S}}^2)$ is smooth.}
\medskip

  The above result will be proved by using the Nash-Moser implicit function theorem; see Section 2.
\medskip

  In order to state our result on asymptotic stability of the radial stationary solution of the problem (1.1), we need to
  introduce some more notations. Let $m$ and $\mu$ be as above. Recall that a bounded domain
  $\Omega\subseteq\mathbb{R}^3$ is called a $C^{m+\mu}$-domain if it is $C^{m+\mu}$-diffeomorphic to the unit sphere $B(0,1)
  \subseteq\mathbb{R}^3$, and $\Omega$ is called a $\dot{C}^{m+\mu}$-domain if it is $\dot{C}^{m+\mu}$-diffeomorphic to the unit
  sphere $B(0,1)$, where $\dot{C}^{m+\mu}$ refers to $m+\mu$-th order little H\"{o}lder continuous class. We use the notation
  $\dot{\mathfrak{D}}^{m+\mu}({\mathbf{R}}^3)$ to denote the Banach manifold of all $\dot{C}^{m+\mu}$-domains in $\mathbb{R}^3$;
  cf. \cite{Cui4} for details. We denote
$$
  \mathfrak{M}:=\dot{\mathfrak{D}}^{m+\mu}({\mathbf{R}}^3), \quad
  \mathfrak{M}_0:=\dot{\mathfrak{D}}^{m+\mu+3}({\mathbf{R}}^3).
$$
  From \cite{Cui4} we know that $\mathfrak{M}_0$ is a $C^3$-embedded Banach submanifold of $\mathfrak{M}$, so that every point in
  $\mathfrak{M}_0$ is $C^3$-differentiable as a point of $\mathfrak{M}$. Given $\Omega\in\mathfrak{M}_0$, the equations
  $(1.1)_1$--$(1.1)_4$ with $\Omega(t)$ replaced by $\Omega$ has a unique solution $(\sigma,p)$ satisfying the following properties:
\begin{equation}
\left\{
\begin{array}{c}
   \sigma,p\in W^{2,q}(\Omega) \;(\forall q\in [1,\infty)), \quad
  \sigma,p\in\dot{C}^{\infty}(\Omega_{\rm liv}\cup\Omega_{\rm nec}),\\
  \sigma\in\dot{C}^{m+3+\mu}(\Omega_{\rm liv}\cup\partial\Omega) \quad \mbox{and} \quad
  p\in\dot{C}^{m+1+\mu}(\Omega_{\rm liv}\cup\partial\Omega),
\end{array}
\right.
\end{equation}
  where $\Omega_{\rm liv}$ and $\Omega_{\rm nec}$ are defined similarly as in (1.3) and (1.4), respectively. We define a vector
  field $\mathscr{F}$ in $\mathfrak{M}$ with domain $\mathfrak{M}_0$ (i.e., $\mathscr{F}:\mathfrak{M}_0\to
  \mathcal{T}_{\mathfrak{M}_0}(\mathfrak{M})$) by letting
\begin{equation}
  \mathscr{F}(\Omega)=-\partial_{\bfn}p|_{\partial\Omega}, \quad  \forall\Omega\in\mathfrak{M}_0.
\end{equation}
  Then the problem (1.1) reduces into the following differential equation in the Banach manifold $\mathfrak{M}$:
\begin{equation}
\left\{
\begin{array}{ll}
   \Omega'(t)=\mathscr{F}(\Omega(t)), &\quad  t>0,\\
    \Omega(0)=\Omega_0. &
\end{array}
\right.
\end{equation}
  We note that when the domain $\Omega(t)$ is determined, the other two components $\sigma,p$ of the solution of the
  problem (1.1) will then follow from solving the elliptic boundary value problem $(1.1)_1$--$(1.1)_4$, and properties
  of $\sigma,p$ are fully determined by those of the domain $\Omega(t)$. Hence, to avoid a very complex statement, we
  only give the precise statement of our result on the initial value problem (1.14), which reads as follows:
\medskip

  {\bf Theorem 1.2}\ \ {\em Let $\mathcal{M}_c$ be the $3$-dimensional submanifold of $\mathfrak{M}_0$ consisting of
  all surface spheres in ${\mathbb{R}}^3$ of radius $R_s$. There exists a constant $\gamma^*>0$ such that if $\gamma
  >\gamma^*$ then the following assertions hold:

  $(1)$\ \ There is a neighborhood $\mathcal{O}$ of $\mathcal{M}_c$ in $\mathfrak{M}_0$ such that for any $\Omega_0\in\mathcal{O}$,
  the initial value problem $(1.14)$ has a unique solution $\Omega\in C([0,\infty),\mathfrak{M}_0)\cap C^1((0,\infty),\mathfrak{M}_0)$.

  $(2)$\ \ There exists a submanifold $\mathcal{M}_s$ of $\mathfrak{M}_0$ of codimension $3$ passing $\Omega_s=B(0,R_s)$
  such that for any $\Omega_0\in\mathcal{M}_s$, the solution of the problem $(1.14)$ satisfies $\displaystyle\lim_{t\to\infty}
  \Omega(t)=\Omega_s$ and, conversely, if the solution of $(1.14)$ satisfies this property then $\Omega_0\in\mathcal{M}_s$.

  $(3)$\ \ For any $\Omega_0\in\mathcal{O}$ there exist unique $x_0\in {\mathbb{R}}^3$ and $\Omega_0'\in\mathcal{M}_s$ such that
  $\Omega_0=x_0+\Omega_0'$ and, for the solution $\Omega=\Omega(t)$ of $(1.14)$, we have
$$
   \lim_{t\to\infty}\Omega(t)=B(x_0,R_s),
$$
  Moreover, convergence rate of the above limit relations is of the form $C\mbox{\rm e}^{-\nu t}$ for some positive constants
  $C$ and $\nu$ depending on $\gamma$. If on the contrary $0<\gamma<\gamma^*$ then the radial stationary solution of the
  problem $(1.14)$ is unstable. $\quad\Box$}
\medskip

  The above result will be proved in Section 4 by using some abstract result for parabolic differential equations in Banach manifold
  recently established in \cite{Cui5}. A crucial step is to show that the representation of the vector field $\mathscr{F}$
  in certain local chart of $\mathfrak{M}$ is smooth, which is ensured by Theorem 1.1.
\medskip

  {\em Remark 1}.\ \ If $R_s>R^*$ then from (1.8) we see that the stationary tumor has a necrotic core, i.e., the set
  $\Omega_{\rm nec}^s={\rm int}\{x\in\Omega_s:\sigma_s(x)=\hat{\sigma}\}$ has a nonempty inner region. By Theorem 1.1 and
  Lemma 2.6 in Section 2 we see that if $\Omega_0$ is sufficiently close to $\Omega_s$ then for any $t>0$ the tumor also
  has a necrotic core, i.e., the set $\Omega_{\rm nec}(t)={\rm int}\{x\in\Omega(t):\sigma(x,t)=\hat{\sigma}\}$ has a nonempty
  inner region. From the assertion (3) and the proof of Theorem 1.1 is is not hard to prove that the following relation
  holds:
$$
   \lim_{t\to\infty}\Omega_{\rm nec}(t)=\Omega_{\rm nec}^s.
$$
  If $R_s<R^*$ then it is also not hard to prove that the tumor does not have a necrotic core for all $t>0$ and the
  results of \cite{CuiEsc2} apply to this situation.
\medskip

  {\em Remark 2}.\ \ As we pointed earlier, properties of the other two components $\sigma,p$ of the solution
  of the problem (1.1) are fully determined by those of the domain $\Omega(t)$; they can be deduced from Theorem 1.2,
  Theorem 1.1 and Lemma 2.6. We omit this discussion here.
\medskip

  The organization of the rest part is as follows. In Section 2 we study the stationary free boundary problem (1.6).
  By using the Nash-Moser implicit function theorem, we prove that the free boundary $\Gamma$ of the problem (1.6)
  is smooth and depends on $\Omega$ smoothly, provided $\Omega$ is a small perturbation of a sphere. We use this
  result to prove the surface tension free part $\pi_0$ of the solution of the second equation in (1.1) has the
  property that the map $\Omega\mapsto \partial_{\bfn}\pi_0|_{\partial\Omega}$ is smooth, even
  if $f$, $g$ are discontinuous functions given by (1.2). In Section 3 we give the proof of Theorem 1.2.

\section{The Proof of Theorem 1.1}
\setcounter{equation}{0}

\hskip 2em
  In this section we give the proof of Theorem 1.1. We shall use Nash-Moser implicit function theorem to prove this theorem.
  Since Nash-Moser implicit function theorem might be not very familiar to the reader, we shall first make a short review to
  this theorem. Hence this section is divided into three subsections: In Subsection 2.1 we make a short review to Nash-Moser
  implicit function theorem. In Subsection 2.2 we give the proof of Theorem 1.1. In Subsection 2.3 we use Theorem 1.1 to prove
  that the representation of the vector field $\mathscr{F}$ in (1.13) in certain local chart of the manifold $\mathfrak{M}$
  is smooth.

\subsection{Review of Nash-Moser implicit function theorem}

\hskip 2em
  Recall that {\em implicit function theorem} in Banach space says that for three Banach spaces $X,Y,Z$ and a $C^1$-mapping
  $F:U\subseteq X\times Y\to Z$, where $U$ is an open subset of $X\times Y$, if $F(x_0,y_0)=0$ for some $(x_0,y_0)\in U$ and
  $\partial_yF(x_0,y_0):Y\to Z$ is an isomorphism of Banach spaces, where $\partial_yF(x,y)$ denotes the partial Fr\'{e}chet
  derivative of $F(x,y)$ in the variable $y$, then there exist $\varepsilon,\delta>0$ with $B(x_0,\varepsilon)\times
  B(y_0,\delta)\subseteq U$ and a unique $C^1$-mapping $f:B(x_0,\varepsilon)\to B(y_0,\delta)$ such that
\begin{equation}
  F(x,f(x))=0
\end{equation}
  for all $x\in B(x_0,\varepsilon)$. This is a fundamental result in the field of mathematical analysis. Naturally one may
  ask whether this result can be extended to general topological vector spaces or more specifically to general Fr\'{e}chet
  spaces. Unfortunately, this is false; cf. counterexamples in Section 5.5 of Part I of \cite{Ham1}. However, if we make suitable
  restrictions to the spaces $X,Y,Z$ as well as the mapping $F$, and suitably strengthen the conditions on $\partial_yF(x,y)$,
  then a generalization named {\em Nash-Moser implicit function theorem} holds. Basic idea of this theorem was first discovered
  by Nash in \cite{Nash} and later formed by Moser in \cite{Moser} into an abstract theorem in functional analysis. An
  excellent exposition of this theorem was given in \cite{Ham1}. In what follows we make a brief review
  to this theorem in the spirit of \cite{Ham1}. Let us start by introducing some basic concepts.

  A {\em graded Fr\'{e}chet space} is a topological vector space $X$ whose topology is defined by an increasing
  sequence of seminorms $\{\|\cdot\|_n\}_{n=0}^{\infty}$, i.e.,
$$
   \|x\|_0\leqslant\|x\|_1\leqslant\|x\|_2\leqslant\cdots\leqslant\|x\|_n\leqslant\cdots, \quad \forall x\in X.
$$
  A linear map $L:X\to Y$ between two graded Fr\'{e}chet spaces $X,Y$ is called a {\em tame linear map} if there exist
  nonnegative integers $n_0,r$ and for each integer $n\geqslant n_0$ a corresponding constant $C_n>0$ such that for any
  integer $n\geqslant n_0$ the following estimate holds:
$$
   \|Lx\|_n\leqslant C_n\|x\|_{n+r}, \quad \forall x\in X.
$$
  A nonlinear map $F:U\subseteq X\to Y$ between two graded Fr\'{e}chet spaces $X,Y$, where $U$ is an open subset of $X$, is
  called a {\em tame nonlinear map} if $F$ is continuous and for any point $x_0\in U$ there exist a neighborhood $U_{x_0}
  \subseteq U$ of $x_0$, nonnegative integers $n_0=n_0(x_0)$, $r=r(x_0)$ and for each integer $n\geqslant n_0$ a corresponding
  constant $C_n(x_0)>0$ such that for any integer $n\geqslant n_0$ the following estimate holds:
$$
   \|F(x)\|_n\leqslant C_n(x_0)(1+\|x\|_{n+r}), \quad \forall x\in U_{x_0}.
$$
  In this case we often simply say that $F$ is tame. If $F:U\subseteq X\to Y$ is smooth and not only itself is tame, but also
  all its derivatives $D^kF$ ($k=1,2,\cdots$) are tame, then we call $F$ a {\em smooth tame map}.

  A graded Fr\'{e}chet space $X$ is called a {\em tame direct summand} of another graded Fr\'{e}chet space $Y$ if there exist
  tame linear maps $F:X\to Y$ and $G:Y\to X$ such that $G(F(x))=x$ for all $x\in X$. For a Banach space $B$, the notation
  $\Sigma(B)$ denotes the space of all sequences $\{x_k\}$ of elements in $B$ such that
$$
   \|\{x_k\}\|_n=\sum_{k=1}^{\infty}e^{nk}\|x_k\|_B<\infty, \quad n=0,1,2,\cdots.
$$
  $\Sigma(B)$ is a graded Fr\'{e}chet space with increasing sequence of seminorms as above. A graded Fr\'{e}chet space $X$ is
  called a {\em tame Fr\'{e}chet space} or simply a {\em tame space} if there exists a Banach space $B$ such that $X$ is a tame
  direct summand of $\Sigma(B)$. In this case we often simply say that $X$ is tame.

  Nash-Moser implicit function theorem mentioned above is concerned with smooth tame maps between tame Fr\'{e}chet spaces.
  The condition for tame Fr\'{e}chet space is quite implicit and hard to verify. Fortunately, we have the following basic results:

  (1) All Banach spaces are tame spaces.

  (2) If $X$ is a compact manifold then $C^{\infty}(X)$ is a tame space.

  (3) If $X$ is a compact manifold with boundary then both $C^{\infty}(X)$ and $C_0^{\infty}(X)$ are tame spaces.

  (4) If $X$ is a compact manifold and $V$ is a vector bundle over $X$ then the space $C^{\infty}(X,V)$ of all smooth sections
  of $V$ over $X$ is a tame space.

  (5) A tame direct summand of a tame space is tame.

  (6) A cartesian product of two tame spaces is tame.

\noindent
  Concerning tame maps and smooth tame maps, we have the following assertions:

  (7) Any continuous map from a graded Fr\'{e}chet space to a Banach space is tame. Any continuous map from a finite
  dimensional space to a graded Fr\'{e}chet space is tame.

  (8) A composition of tame maps is tame.

  (9) Let $X$ be a compact manifold and $V,W$ be vector bundles over $X$. Let $U$ be an open subset of $V$ and $p:U\subseteq
  V\to W$ be a smooth map of $U$ into $W$ which takes fibres into fibres. Let $\tilde{U}\subseteq C^{\infty}(X,V)$ be the
  set of smooth sections of $V$ over $X$ whose image lies in $U$. Then $\tilde{U}$ is an open subset of $C^{\infty}(X,V)$
  and the map $P:\tilde{U}\subseteq C^{\infty}(X,V)\to C^{\infty}(X,W)$ defined by $Pf(x)=p(f(x))$ (for $f\in\tilde{U}$),
  called {\em nonlinear vector bundle operator}, is tame.

  (10) Let $X,V,W$ be as above, $m$ a positive integer and $U$ an open subset of $C^{\infty}(X,V)$. A smooth {\em nonlinear
  partial differential operator} $P$ of order $m$ from $V$ to $W$ is a map $P:U\subseteq C^{\infty}(X,V)\to C^{\infty}(X,W)$
  such that for any $f\in U$ and $x\in X$, $Pf(x)$ is a smooth function of $f(x)$ and partial derivatives of $f$ at $x$ of
  degree at most $m$ in any local charts. A smooth nonlinear partial differential operator is a smooth tame map.
\medskip

  For proofs of the above assertions, we refer the reader to see \cite{Ham1}.

  Nash-Moser implicit function theorem reads as follows:
\medskip

  {\bf Theorem 2.1} \ {\em Let $X,Y,Z$ be tame Fr\'{e}chet spaces and $F:U\subseteq X\times Y\to Z$ a smooth tame map,
  where $U$ is an open subset of $X\times Y$. Let $(x_0,y_0)\in U$ be such that $F(x_0,y_0)=0$. Assume that there exists
  a smooth tame map $A:U(\subseteq X\times Y)\times Z\to Y$ of the form $A(x,y,z)=L(x,y)z$, where for each $(x,y)\in U$,
  $L(x,y)$ is a linear map from $Z$ to $Y$, such that
\begin{equation}
   L(x,y)\partial_yF(x,y)u=u \quad \mbox{and} \quad \partial_yF(x,y)L(x,y)z=z
\end{equation}
  for all $(x,y)\in U$, $u\in Y$ and $z\in Z$. Then there exist neighborhoods $B_1, B_2$ of $x_0$ and $y_0$, respectively,
  with the property that $B_1\times B_2\subseteq U$, and a smooth tame map $f:B_1\to B_2$, such that
$$
   f(x_0)=y_0  \quad \mbox{and} \quad F(x,f(x))=0
$$
  for all $x\in B_1$. Moreover, for any $x\in B_1$, $y=f(x)$ is the unique solution of the equation $F(x,y)=0$ in $B_2$.}
\medskip

  For the proof of the above theorem, we refer the reader to see Theorems 3.3.1 and 3.3.3 in Part III of \cite{Ham1}. We
  note that in those theorems some quadratic error terms are included in the right-hand sides of the two equations in (2.2).
  Here we take such terms to be identically vanishing.
\medskip

   {\em Remark}.  Comparing Nash-Moser implicit function theorem with the implicit function theorem in Banach space, we see
   that a significant difference between these theorems is that in the Fr\'{e}chet space case, the partial derivative
   $\partial_yF(x,y)$ of $F(x,y)$ should be invertible not merely at the single point $(x_0,y_0)$ as in the Banach space case,
   but at all points in a neighborhood of this point. Partial reason for this difference to occur is due to the fact that for
   two Banach spaces $X$ and $Y$, the set of invertible continuous linear maps from $X$ to $Y$ is an open subset of $L(X,Y)$,
   whereas if $X$ and $Y$ are Fr\'{e}chet spaces, this is not the case, even if they are tame; cf. the counterexample 5.3.3
   in Part I of \cite{Ham1}.

\subsection{The Proof of Theorem 1.1}

\hskip 2em
  In this subsection we give the proof of Theorem 1.1.

  Let $R>R^*$ be given and set $K=K(R)$. Let $m,\mu$ be as in Theorem 1.1, i.e., $m$ is an integer not less than $2$ and
  $0<\mu<1$. We know that $C^{\infty}({\mathbf{S}}^2)$ with the family of seminorms $\{\|\cdot\|_{C^({\mathbf{S}}^2)}\}\cup
  \{\|\cdot\|_{C^{k+\mu}({\mathbf{S}}^2)}\}_{k=1}^{\infty}$ is a tame Frech\'{e}t space. We also regard the Banach space
  $C^{m+\mu}({\mathbf{S}}^2)$ as a tame Frech\'{e}t space. For sufficiently small $\delta,\delta'>0$ we denote
$$
  O_{\delta}=\{\rho\in C^{m+\mu}({\mathbf{S}}^2):\|\rho\|_{C^{m+\mu}({\mathbf{S}}^2)}<\delta\}, \quad
  O_{\delta'}'=\{\eta\in C^{\infty}({\mathbf{S}}^2):\|\eta\|_{C^{m+\mu}({\mathbf{S}}^2)}<\delta'\};
$$
  they are open subsets of $C^{m+\mu}({\mathbf{S}}^2)$ and $C^{\infty}({\mathbf{S}}^2)$, respectively. We define a map
  $A: O_{\delta}\times O_{\delta'}'\subseteq C^{m+\mu}({\mathbf{S}}^2)\times C^{\infty}({\mathbf{S}}^2)\to
  C^{\infty}({\mathbf{S}}^2)$ as follows: Given $\rho\in O_{\delta}$ and $\eta\in O_{\delta'}'$, let $\sigma=
  \sigma(r,\omega;\rho,\eta)$ be the unique solution of the equations $(1.11)_1$, $(1.11)_2$ and $(1.11)_4$, and define
\begin{equation}
  A(\rho,\eta)=[\omega\mapsto\sigma(K[1+\eta(\omega)],\omega;\rho,\eta)-\hat{\sigma},\omega\in{\mathbf{S}}^2].
\end{equation}
  Clearly $A(0,0)=0$. We shall prove that if $\delta,\delta'$ are sufficiently small then there exists a unique
  smooth mapping $\varphi:O_{\delta}\to O_{\delta'}'$ such that $A(\rho,\varphi(\rho))=0$ for all $\rho\in O_{\delta}$.
\medskip

  {\bf Lemma 2.2}\ \ {\em $A$ is a smooth tame map.}
\medskip

  {\bf Proof}.\ \ Choose a function $\phi\in C^{\infty}[ K,R]$ such that it satisfies the following conditions:
$$
  0\leqslant\phi\leqslant 1; \quad \phi(R)=\phi( K)=1; \quad
  \phi(t)=0 \;\; \mbox{for}\;\,\frac{3}{4} K+\frac{1}{4}R\leqslant t\leqslant\frac{1}{4} K+\frac{3}{4}R;
$$
$$
  \phi'(t)\leqslant 0 \;\; \mbox{for}\;\, K\leqslant t\leqslant\frac{3}{4} K+\frac{1}{4}R; \quad
  \phi'(t)\geqslant 0 \;\; \mbox{for}\;\,\frac{1}{4} K+\frac{3}{4}R\leqslant t\leqslant R.
$$
  Let $M_0=\displaystyle\max_{ K\leqslant t\leqslant R}|\phi'(t)|$ and assume $\delta,\delta'$ are small enough such
  that $\delta<(1+M_0R)^{-1}$, $\delta'<(1+M_0 K)^{-1}$ and $\max\{\delta,\delta'\}<\displaystyle\frac{1}{3}
  \frac{R- K}{R+ K}$. Consider the variable transformation $y=\Psi_{\rho,\eta}(x)$ from $\overline{D}_{\rho,\eta}$ to
  $\overline{D}$, where for $x\in\overline{D}_{\rho,\eta}$,
\begin{equation}
  \Psi_{\rho,\eta}(x)=\left\{
\begin{array}{ll}
   \displaystyle x-R\rho(\omega)\phi\Big(\frac{r}{1+\rho(\omega)}\Big)\omega &\quad\;\;\mbox{if} \;\; r\geqslant
   \frac{1}{2}( K+R),
\\ [0.3cm]
   \displaystyle x- K\eta(\omega)\phi\Big(\frac{r}{1+\eta(\omega)}\Big)\omega &\quad\;\;\mbox{if} \;\; r<\frac{1}{2}( K+R).
\end{array}
\right.
\end{equation}
  It is easy to see that $\Psi_{\rho,\eta}$ is a $C^{m+\mu}$ diffeomorphism from $\overline{D}_{\rho,\eta}$ onto
  $\overline{D}$, and
$$
  \Psi_{\rho,\eta}(x)=x \quad \mbox{for}\;\;\frac{1}{6}\frac{(3K\!+\!R)(K\!+\!2R)}{K\!+\!R}\leqslant
  |x|\leqslant\frac{1}{6}\frac{(K\!+\!3R)(2K\!+\!R)}{K\!+\!R}.
$$
  Moreover, denoting
$$
  E_{\eta}=\{x\in{\mathbf{R}}^3:  K[1+\eta(\omega)]<r<\frac{1}{2}( K+R)\}, \quad
  E=\{x\in{\mathbf{R}}^3:  K<r<\frac{1}{2}(K+R)\},
$$
  we see that the restriction of $\Psi_{\rho,\eta}$ on $\overline{E}_{\eta}$ is a $C^{\infty}$-diffeomorphism from
  $\overline{E}_{\eta}$ onto $\overline{E}$, due to the facts that $\eta\in C^{\infty}({\mathbf{S}}^2)$ and that this
  restriction is independent of $\rho$. Because of the latter property, we re-denote the restriction of $\Psi_{\rho,\eta}$
  on $\overline{E}_{\eta}$ as $\Psi_{\eta}$, and denote by $\psi_{\eta}$ the restriction of $\Psi_{\eta}$ to $\Gamma_{\eta}$,
  which is clearly a $C^{\infty}$-diffeomorphism from $\Gamma_{\eta}$ onto $\Gamma_0$. Now define operators
  $\mathscr{A}(\rho,\eta):C^{m+\mu}(\overline{D})\cap C^{\infty}(\overline{E})\to C^{m-2+\mu}(\overline{D})\cap
  C^{\infty}(\overline{E})$, $\mathscr{A}(\eta):C^{\infty}(\overline{E})\to C^{\infty}(\overline{E})$ and
  $\mathscr{N}(\eta):C^{\infty}(\overline{E})\to C^{\infty}(\Gamma_0)$ respectively as follows:
$$
  \mathscr{A}(\rho,\eta)u=[\Delta(u\circ\Psi_{\rho,\eta})]\circ\Psi_{\rho,\eta}^{-1} \quad\;\;
   \mbox{for}\;\,u\in C^{m+\mu}(\overline{D})\cap C^{\infty}(\overline{E}),
$$
$$
  \mathscr{A}(\eta)u=[\Delta(u\circ\Psi_{\eta})]\circ\Psi_{\eta}^{-1} \quad\;\;
   \mbox{for}\;\,u\in C^{\infty}(\overline{E}),
$$
$$
  \mathscr{N}(\eta)u=[\partial_r(u\circ\Psi_{\eta})|_{\Gamma_{\eta}}]\circ\psi_{\eta}^{-1} \quad\;\;
  \mbox{for}\;\,u\in C^{\infty}(\overline{E}).
$$
  Let $u=\sigma\circ\Psi_{\rho,\eta}^{-1}$. After the variable transformation $x\mapsto\Psi_{\rho,\eta}(x)$, the problem
  $(1.11)_1$, $(1.11)_2$ and $(1.11)_4$ transforms into the following problem:
\begin{equation}
\left\{
\begin{array}{rll}
   \mathscr{A}(\rho,\eta)u=&u &  \quad  \mbox{in}\;\; D,\\
   u=&1 &  \quad  \mbox{on}\;\; S_0,\\
   \mathscr{N}(\eta)u=&0 & \quad  \mbox{on}\;\; \Gamma_0.
\end{array}
\right.
\end{equation}
  It is clear that all coefficients of the operator $\mathscr{A}(\rho,\eta)$ belong to $C^{m-2+\mu}(\overline{D})$ and the
  $C^{m-2+\mu}(\overline{D})$-norms of all these coefficients are bounded with a constant depending only on $m,\mu,\delta,\delta'$,
  and similarly all coefficients of the operator $\mathscr{N}(\eta)$ belong to $C^{m-1+\mu}(\Gamma_0)$ and the
  $C^{m-1+\mu}(\Gamma_0)$-norms of all these coefficients are also bounded with a constant depending only on $m,\mu,\delta,\delta'$.
  Besides, it is also clear that the smallest eigenvalue of the second-order coefficient matrix of the operator $-\mathscr{A}(\rho,\eta)$
  is bounded below by a positive constant depending only on $m,\mu,\delta,\delta'$. Hence, by a standard result in the theory of
  elliptic boundary value problems we see that the solution $u$ satisfies the following estimates:
\begin{equation}
   0<u\leqslant 1; \qquad
   \|u\|_{C^{m+\mu}(\overline{D})}\leqslant C(m,\mu,\delta,\delta')<\infty.
\end{equation}
  Next, choose a number $\varepsilon>0$ sufficiently small such that $\varepsilon<\displaystyle\frac{K}{18}\frac{R\!-\!K}{R\!+\!K}$
  and let
$$
  \Lambda_j=\{x\in{\mathbf{R}}^3:  \frac{1}{2}( K+R)-j\varepsilon<r<\frac{1}{2}( K+R)+j\varepsilon\}, \quad
  j=1,2,3.
$$
  Note that $\Lambda_1\subset\subset\Lambda_2\subset\subset\Lambda_3$.
  Since in $\Psi_{\rho,\eta}(x)=x$ for $x\in\Lambda_3$, from the equation $\mathscr{A}(\rho,\eta)u=u$ we have
$$
  \Delta u=u \quad \mbox{in} \;\; \Lambda_3.
$$
  From this fact and the first estimate in (2.6), by using some standard estimates for elliptic equations we get:
\begin{equation}
   \|u\|_{C^{k+\mu}(\overline{\Lambda}_2)}\leqslant C(k), \quad k=m,m+1,m+2,\cdots.
\end{equation}
  In $E$ we have
$$
\left\{
\begin{array}{rll}
   \mathscr{A}(\eta)u=&u &  \quad  \mbox{in}\;\; E,\\
   \mathscr{N}(\eta)u=&0 & \quad  \mbox{on}\;\; \Gamma_0.
\end{array}
\right.
$$
  The operators $\mathscr{A}(\eta)$ and $\mathscr{N}(\eta)$ respectively have the following forms:
$$
   \mathscr{A}(\eta)u=\sum_{i,j=1}^3a_{ij}(x,\eta,\nabla\eta)\partial_{ij}^2u+
   \sum_{i=1}^3b_{i}(x,\eta,\nabla\eta,\nabla^2\eta)\partial_{i}u,
$$
$$
   \mathscr{N}(\eta)u=\sum_{i=1}^3c_{i}(x,\eta,\nabla\eta)\partial_{i}u,
$$
  where $a_{ij}(x,\eta,\nabla\eta)$'s are quadratic functions in $\nabla\eta$ with coefficients being smooth functions of
  $x$ and $\eta$, $b_{i}(x,\eta,\nabla\eta,\nabla^2\eta)$'s are sums of linear functions in $\nabla^2\eta$ and quadratic
  functions in $\nabla\eta$ with coefficients being smooth functions of $x$ and $\eta$, and $c_{i}(x,\eta,\nabla\eta)$'s
  are linear functions in $\nabla\eta$ with coefficients being smooth functions of $x$ and $\eta$. Using these facts, the
  estimates in (2.6) and some standard arguments as in the proofs of higher-order interior and boundary regularity estimates
  for elliptic boundary value problems, we see that the following estimates hold:
\begin{equation}
   \|u\|_{C^{k+\mu}(\overline{E\backslash\Lambda}_1)}\leqslant C(k,\mu,\delta,\delta')(1+\|\eta\|_{C^{k+\mu}({\mathbf{S}}^2)}),
   \quad k=m+1,m+2,\cdots.
\end{equation}
  Combining the estimates (2.6), (2.7) and (2.8), we get the following estimates:
\begin{equation}
  \|u\|_{C^{m+\mu}(\overline{D})}+\|u\|_{C^{k+\mu}(\overline{E})}\leqslant C_k(1+\|\eta\|_{C^{k+\mu}({\mathbf{S}}^2)}),
  \quad k=m,m+1,m+2,\cdots.
\end{equation}
  Next, for any positive integer $p$, if we denote by $u_p$ any one of the partial Fr\'{e}chet derivatives of order $p$ of
  $u$ in $\rho,\eta$ (recall that $u_p$ is a continuous $p$-linear operator in $C^{\infty}({\mathbf{S}}^2)$ with value in
  $C^{m+\mu}(\overline{D})\cap C^{\infty}(\overline{E})$), then from (2.5) we easily see that for any
  $\xi_1,\xi_2,\cdots,\xi_p\in C^{\infty}({\mathbf{S}}^2)$, $u_p(\xi_1,\xi_2,\cdots,\xi_p)$ satisfies an elliptic boundary
  value problem of the following form:
$$
\left\{
\begin{array}{rll}
   \mathscr{A}(\rho,\eta)u_p(\xi_1,\xi_2,\cdots,\xi_p)=&u_p(\xi_1,\xi_2,\cdots,\xi_p)
   +\mathcal{F}_p(\rho,\eta,u,Du_1,\cdots,D^{p-1}u,\xi_1,\xi_2,\cdots,\xi_p) &  \quad  \mbox{in}\;\; D,\\
   u_p(\xi_1,\xi_2,\cdots,\xi_p)=&0 &  \quad  \mbox{on}\;\; S_0,\\
   \mathscr{N}(\eta)u_p(\xi_1,\xi_2,\cdots,\xi_p)=&
   \mathcal{G}_p(\rho,\eta,u,Du_1,\cdots,D^{p-1}u,\xi_1,\xi_2,\cdots,\xi_p) & \quad  \mbox{on}\;\; \Gamma_0,
\end{array}
\right.
$$
  where $D^iu$ denotes the set of partial Fr\'{e}chet derivatives of $u$ in $\rho,\eta$ of order $i$, $i=1,2,\cdots,p-1$,
  $\mathcal{F}_p(\rho,\eta,u,Du_1,\cdots,D^{p-1}u,\xi_1,\xi_2,\cdots,\xi_p)$ and $\mathcal{G}_p(\rho,\eta,u,Du_1,\cdots,
  D^{p-1}u,\xi_1,\xi_2,\cdots,\xi_p)$ denotes functionals in $\rho,\eta,u,Du_1,\cdots,D^{p-1}u,\xi_1,\xi_2,\cdots,\xi_p$,
  both linear in $(u,Du_1,\cdots,D^{p-1}u)$ and $p$-linear in\\
  $(\xi_1,\xi_2,\cdots,\xi_p)$. Using these facts, and a similar argument as
  above and induction in $p$, we see that similar estimates as in (2.9) also hold for $u_p$ for any $p\geqslant 1$, i.e.,
  for any $(\rho,\eta)\in O_{\delta}\times O_{\delta'}'$ and any $\xi_1,\xi_2,\cdots,\xi_p\in C^{\infty}({\mathbf{S}}^2)$,
\begin{align}
  &\|u_p(\xi_1,\xi_2,\cdots,\xi_p)\|_{C^{m+\mu}(\overline{D})}+\|u_p(\xi_1,\xi_2,\cdots,\xi_p)\|_{C^{k+\mu}(\overline{D})}
\nonumber\\
  \leqslant&\displaystyle C_{p,k}(1+\|\eta\|_{C^{k+\mu}({\mathbf{S}}^2)})\sum_{i=1}^p\|\xi_1\|_{C^0({\mathbf{S}}^2)}\cdots
  \|\xi_{i-1}\|_{C^0({\mathbf{S}}^2)}\|\xi_{i}\|_{C^{2}({\mathbf{S}}^2)}\|\xi_{i+1}\|_{C^0({\mathbf{S}}^2)}\cdots
  \|\xi_p\|_{C^0({\mathbf{S}}^2)}
\nonumber\\
  &+C_{p,k}'\sum_{i=1}^p\|\xi_1\|_{C^0({\mathbf{S}}^2)}\cdots
  \|\xi_{i-1}\|_{C^0({\mathbf{S}}^2)}\|\xi_{i}\|_{C^{k+\mu}({\mathbf{S}}^2)}\|\xi_{i+1}\|_{C^0({\mathbf{S}}^2)}\cdots
  \|\xi_p\|_{C^0({\mathbf{S}}^2)},
\nonumber\\
  &\qquad\qquad\qquad\qquad\qquad\qquad  k=m,m+1,m+2,\cdots.
\end{align}
  From (2.9) and (2.10) we see that the solution map $(\rho,\eta)\mapsto u$ of the problem (2.8) is a smooth tame map. Since
$$
  A(\rho,\eta)=\Big[\omega\mapsto u(K\omega)-\hat{\sigma},
  \omega\in{\mathbf{S}}^2\Big],
$$
  the desired assertion immediately follows. $\quad\Box$
\medskip

  {\bf Lemma 2.3}\ \ {\em There exist $\delta,\delta'>0$ sufficiently small such that for all $\rho\in O_{\delta}$
  and $\eta\in O_{\delta'}'$, $\partial_{\eta}A(\rho,\eta)$ is invertible, and the map $(\rho,\eta,\xi)\mapsto
  \partial_{\eta}A(\rho,\eta)^{-1}\xi$ is a smooth tame map.}
\medskip

  {\bf Proof}.\ \ We note that from the definition of $A(\rho,\eta)$ (see (2.3)), we see that if as in (2.3) we denote
  by $\sigma=\sigma(r,\omega;\rho,\eta)$ the solution of the problem $(1.11)_1$, $(1.11)_2$ and $(1.11)_4$, then for
  any $\xi\in C^{\infty}({\mathbf{S}}^2)$ we have
$$
\begin{array}{rl}
   \partial_{\eta}A(\rho,\eta)\xi=&\partial_{\eta}\sigma(K[1+\eta(\omega)],\omega;\rho,\eta)\xi
  +\partial_r\sigma(K[1+\eta(\omega)],\omega;\rho,\eta)K\xi
\\ [1mm]
   =&\partial_{\eta}\sigma(K[1+\eta(\omega)],\omega;\rho,\eta)\xi \quad (\mbox{by $(1.11)_4$}).
\end{array}
$$
  Let $v(r,\omega)=\partial_{\eta}\sigma(r,\omega;\rho,\eta)\xi$. Then $\partial_{\eta}A(\rho,\eta)\xi=[\omega\mapsto
  v(K[1+\eta(\omega)],\omega),\omega\in{\mathbf{S}}^2]$. A simple computation shows that $v$ is the solution of the
  following problem:
$$
  \Delta v=v \;\;  \mbox{in} \;\; D_{\rho,\eta}, \quad  v=0 \;\;  \mbox{on}\;\; S_{\rho},  \quad
   \partial_rv=-K\partial_r^2\sigma|_{\Gamma_{\eta}}\xi\;\;\;  \mbox{on} \;\; \Gamma_{\eta}.
$$
  If $\rho=0$ and $\eta=0$ then $\sigma=U(r,R)$, so that $\partial_r^2\sigma(\cdot;0,0)|_{\Gamma_{\eta}}=\partial_r^2U(K,R)=U(K,R)=\hat{\sigma}$.
  We now choose $\delta,\delta'>0$ sufficiently small such that for all $\rho\in O_{\delta}$ and $\eta\in O_{\delta'}'$ there hold $(1/2)\hat{\sigma}
  \leqslant\partial_r^2\sigma(\cdot;\rho,\eta)|_{\Gamma_{\eta}}\leqslant 2\hat{\sigma}$. Then for any $\rho\in O_{\delta}$ and $\eta\in O_{\delta'}'$
  the operator $\partial_{\eta}A(\rho,\eta)$ is invertible, with
$$
  [\partial_{\eta}A(\rho,\eta)]^{-1}\zeta=-\frac{\partial_rw(\cdot;\rho,\eta,\zeta)}{K\partial_r^2\sigma(\cdot;\rho,\eta)}\Big|_{\Gamma_{\eta}},
  \quad \forall\zeta\in C^{\infty}({\mathbf{S}}^{n-1}),
$$
  where $w=w(\cdot;\rho,\eta,\zeta)$ is the solution of the following problem:
$$
  \Delta w=w \;\;  \mbox{in} \;\; D_{\rho,\eta}, \quad  w=0 \;\;  \mbox{on}\;\; S_{\rho},  \quad
  w=\zeta\;\;  \mbox{on} \;\; \Gamma_{\eta}.
$$
  Here $w=\zeta$ on $\Gamma_{\eta}$ means that $w|_{r=K(1+\eta(\omega))}=\zeta(\omega)$ for all $\omega\in {\mathbf{S}}^2$.
  Similarly as before we can prove the map $(\rho,\eta,\zeta)\mapsto\partial_rw(\cdot;\rho,\eta,\zeta)$ is a smooth
  tame map. Since we have known that the map $(\rho,\eta)\mapsto\sigma(\cdot;\rho,\eta)$ is a smooth tame
  map, it follows that the map $(\rho,\eta)\mapsto\partial_r^2\sigma(\cdot;\rho,\eta)$ is also a smooth tame map
  (by Theorem 2.2.6 in Part II of \cite{Ham1}), and, consequently, the map $(\rho,\eta,\zeta)\mapsto
  \partial_rw(\cdot;\rho,\eta,\zeta)/\partial_r^2\sigma(\cdot;\rho,\eta)$ is a smooth tame map. Hence the desired assertion
  follows. $\quad\Box$
\medskip

  {\bf Proof of Theorem 1.1}.\ \ Having proved Lemmas 2.2 and 2.3, by using the Nash-Moser implicit function theorem
  we conclude that by choosing $\delta,\delta'>0$ smaller when necessary, it follows that for any $\rho\in O_{\delta}$
  there exists $\eta\in O_{\delta'}'$ such that it is the unique solution of the equation $A(\rho,\eta)=0$ in $O_{\delta'}'$,
  and the map $\rho\mapsto\eta$ from $O_{\delta}\subseteq C^{m+\mu}({\mathbf{S}}^2)$ to $O_{\delta'}'\subseteq
  C^{\infty}({\mathbf{S}}^2)$ is a smooth tame map. This shows that for any $\rho\in C^{m+\mu}({\mathbf{S}}^2)$ with
  $\|\rho\|_{C^{m+\mu}({\mathbf{S}}^2)}<\delta$, the free boundary $\Gamma_{\eta}$ is smooth and the mapping $\rho\mapsto\eta$
  is also smooth. Having proved smoothness of the free boundary $\Gamma_{\eta}$, the assertion $\sigma\in
  C^{m+\mu}(\overline{D}_{\rho,\eta})\cap C^{\infty}(\overline{D}_{\rho,\eta}\backslash S_{\rho})$ follows immediately.
  This proves Theorem 1.1. $\quad\Box$
\medskip

  {\em Remark}.\ \ The reader might argue why we don't use the implicit function theorem in Banach space to prove
  theorem 1.1, as in the proof of Lemma 5.1 of \cite{Cui3}. The reason is that the method used in the proof of
  Lemma 5.1 of \cite{Cui3} does not work to the present problem, due to the fact that $\sigma'(K)=0$, where
  $\sigma(r)=U(r,R)$.
\medskip

\subsection{Smoothness of a map related to the free boundary problem (1.1)}

\hskip 2em
  In this subsection we study the obstacle problem (1.6). The purpose is through a such study to prove that for the solution
  $(\sigma,\pi_0)$ of the boundary value problem
\begin{equation}
\left\{
\begin{array}{rll}
   \Delta\sigma=&f(\sigma) &\quad\;\; \mbox{in} \;\; \Omega,\\
   -\Delta\pi_0=&g(\sigma)   &\quad\;\; \mbox{in} \;\; \Omega,\\
   \sigma=&1   &\quad\;\; \mbox{on} \;\; \partial\Omega,\\
   \pi_0=&0  &\quad\;\; \mbox{on} \;\; \partial\Omega,
\end{array}
\right.
\end{equation}
  where $f$, $g$ are the discontinuous functions given in (1.2) (with $\lambda=1$), the mapping $\Omega\mapsto
  \partial_{\bfn}\pi_0|_{\partial\Omega}$ from a neighborhood of a sphere in $\mathfrak{M}_0:=
  \dot{\mathfrak{D}}^{m+3+\mu}({\mathbf{R}}^3)\subseteq\mathfrak{M}:=\dot{\mathfrak{D}}^{m+\mu}({\mathbf{R}}^3)$ to
  $\mathcal{T}_{\mathfrak{M}_0}(\mathfrak{M})$ is smooth, i.e., representation of this mapping in some regular local chart
  of $\mathfrak{M}$ at every sphere is smooth, where $\partial_{\bfn}$ denotes the derivative in the outward normal direction
  $\bfn$ of $\partial\Omega$. This result is crucial in the study of the free boundary problem (1.1) to be given in the next
  section. We note that since the functions $f$, $g$ are discontinuous, such a result apparently looks unbelievable.

  We point out that although here we only consider the three dimension case, a similar discussion also works for general
  dimension $n\geqslant 2$ case; in order to do so the discussion in \cite{Cui2} must be first extended, which is not hard.

  {\bf Lemma 2.4}\ \ {\em Let $m\in\mathbb{N}$, $m\geqslant$, and $0<\mu<1$ be given. For $\rho\in C^{m+\mu}({\mathbf{S}}^2)$
  with $\|\rho\|_{C^{m+\mu}({\mathbf{S}}^2)}$ sufficiently small, the problem $(1.6)$ with $\Omega=\Omega_{\rho}$ has a unique
  solution $\sigma\in\displaystyle\bigcap_{1\leqslant q<\infty}W^{2,q}(\Omega)$, and the solution has the following properties:

  $(1)$\ $\hat{\sigma}\leqslant\sigma(x)\leqslant1$ for all $x\in\overline{\Omega}$.

  $(2)$\ There exists $\eta\in C^{\infty}({\mathbf{S}}^2)$ with $\|\eta\|_{C^{2+\mu}({\mathbf{S}}^2)}$ small, such that
  $\sigma(x)=\hat{\sigma}$ for $x\in\overline{\Omega}_{\nec}$, where $\Omega_{\nec}=\{x\in\mathbb{R}^3:r<K[1+\eta(\omega)]\}$,
  and $\hat{\sigma}<\sigma(x)\leqslant1$ for $x\in\Omega_{\liv}\cup\partial\Omega$, where $\Omega_{\liv}=\Omega\backslash
  \overline{\Omega}_{\nec}$.

  $(3)$\ The map $\rho\mapsto\eta$ from a small neighborhood of the origin of the Banach space $C^{m+\mu}({\mathbf{S}}^2)$
  to the Frech\'{e}t space $C^{\infty}({\mathbf{S}}^2)$ is smooth.

  $(4)$\ $\sigma|_{\overline{\Omega}_{\liv}}\in C^{\infty}(\Omega_{\liv}\cup\Gamma)\cup C^{m+\mu}(\overline{\Omega}_{\liv})$,
  where $\Gamma=\partial\Omega_{\nec}$, and the map $\rho\mapsto\sigma|_{\overline{\Omega}_{\liv}}$ is smooth in the following
  sense: Let $\Psi_{\rho,\eta}$ be the Hanzawa transformation given by $(2.4)$. Then the map $\rho\mapsto\sigma\circ
  \Psi_{\rho,\eta}^{-1}$ from a small neighborhood of the origin of the Banach space $C^{m+\mu}({\mathbf{S}}^2)$ to the
  Frech\'{e}t space $C^{\infty}(\Omega_{\liv}\cup\Gamma)\cup C^{m+\mu}(\overline{\Omega}_{\liv})$ is smooth. Note that
  $\Psi_{\rho,\eta}$ depends on $\rho$ smoothly.}
\medskip

  {\em Proof}.\ \ We divide the proof into three steps.

  Step 1:\ By Theorem 1.1 we know that there exists $\delta>0$ sufficiently small, such that given $\rho\in
  C^{m+\mu}({\mathbf{S}}^2)$ $(0<\mu<1)$ with $\|\rho\|_{C^{m+\mu}({\mathbf{S}}^2)}<\delta$, the problem $(1.11)$ has a unique
  solution $(\sigma,\eta)$ with $\eta\in C^{\infty}({\mathbf{S}}^2)$ and $\sigma\in C^{2+\mu}(\overline{D}_{\rho,\eta})\cap
  C^{\infty}(\overline{D}_{\rho,\eta}\backslash S_{\rho})$, and the mapping $\rho\mapsto\eta$ from the open set
  $\|\rho\|_{C^{m+\mu}({\mathbf{S}}^2)}<\delta$ in $C^{m+\mu}({\mathbf{S}}^2)$ to the Frech\'{e}t space
  $C^{\infty}({\mathbf{S}}^2)$ is smooth. Moreover, the proof of that theorem also ensures that the map $\rho\mapsto\sigma\circ
  \Psi_{\rho,\eta}^{-1}$ from the open set $\|\rho\|_{C^{m+\mu}({\mathbf{S}}^2)}<\delta$ in $C^{m+\mu}({\mathbf{S}}^2)$ to
  $C^{m+\mu}(\overline{D})\cap C^{\infty}(\overline{D}\backslash S_0)$ is smooth. Note that by maximum principle it is clear
  that $\hat{\sigma}\leqslant\sigma(x)\leqslant1$ for all $x\in\overline{D}_{\rho,\eta}$.

  Step 2:\ Let $\nu$ be the outward unit normal field of the inner part boundary $\Gamma_{\eta}$ of $D_{\rho,\eta}$.
  We prove that under the assumption that $\rho,\eta\in C^{m+\mu}({\mathbf{S}}^2)$ and $\|\eta\|_{C^{m+\mu}({\mathbf{S}}^2)}$
  is sufficiently small, the problem (1.11) is equivalent to the following problem:
\begin{equation}
\left\{
\begin{array}{rll}
   \Delta\sigma=&\sigma &  \quad  \mbox{in}\;\; D_{\rho,\eta},\\
   \sigma=&1 &  \quad  \mbox{on}\;\; S_{\rho},\\
   \sigma=&\hat{\sigma} &  \quad  \mbox{on}\;\; \Gamma_{\eta},\\
    \partial_{\nu}\sigma=&0 & \quad  \mbox{on}\;\; \Gamma_{\eta}.
\end{array}
\right.
\end{equation}
  Indeed, the condition $\sigma=\hat{\sigma}$ on $\Gamma_{\eta}$ implies $\nabla_{\omega}\sigma(\eta(\omega),\omega)=
  -\partial_r\sigma(\eta(\omega),\omega)\nabla_{\omega}\eta(\omega)$ for $\omega\in{\mathbf{S}}^{n-1}$. Since
  $\nu=[\omega-(1/r)\nabla_{\omega}\eta(\omega)]/\sqrt{1+(1/r^2)|\nabla_{\omega}\eta(\omega)|^2}|_{r=\eta(\omega)}$ and
  $\nabla\sigma=(\partial_{r}\sigma)\omega+(1/r)\nabla_{\omega}\sigma$, so that
$$
  \partial_{\nu}\sigma|_{r=\eta(\omega)}=\frac{\partial_{r}\sigma-(1/r^2)\nabla_{\omega}\eta(\omega)
  \cdot\nabla_{\omega}\sigma}{\sqrt{1+(1/r^2)|\nabla_{\omega}\eta(\omega)|^2}}\Big|_{r=\eta(\omega)},
$$
  we see that the condition $\sigma=\hat{\sigma}$ on $\Gamma_{\eta}$ implies $\partial_{\nu}\sigma|_{r=\eta(\omega)}
  =\sqrt{1+(1/r^2)|\nabla_{\omega}\eta(\omega)|^2}\partial_{r}\sigma|_{r=\eta(\omega)}$. Hence the problems (1.11) and
  (2.12) are equivalent.

  Step 3:\ It is easy to see that if $\sigma$ is a solution of (2.12) then by extending it into the whole domain
  $\overline{\Omega}=\overline{\Omega}_{\rho}$ such that it identically takes the value $\hat{\sigma}$ in $\Omega_{\nec}=
  \Omega\backslash D_{\rho,\eta}$, then after such extension $\sigma\in\displaystyle\bigcap_{1\leqslant q<\infty}
  W^{2,q}(\Omega)$ and it  is a solution of (1.6). By the above two steps, it follows that given $\rho\in
  C^{m+\mu}({\mathbf{S}}^2)$ with $\|\rho\|_{C^{m+\mu}({\mathbf{S}}^2)}<\delta$, the problem (1.6) has a unique
  solution $\sigma\in\displaystyle\bigcap_{1\leqslant q<\infty}W^{2,q}(\Omega)$. Clearly, this solution possesses the
  properties (1)--(4). Since the $f$ is a monotone nondecreasing function, by using the weak maximum principle
  (cf. Theorem 2.3 in \cite{Tro}) we see that the solution of the problem (1.6) is unique. This completes the proof of
  Lemma 2.4. $\quad\Box$
\medskip

  {\em Remark}.\ As an immediate corollary of the above lemma, we see that as far as strong solution is concerned, the
  problems (1.6) and (1.11) are equivalent.
\medskip

  Let us now consider the problem (2.11). Let $\delta$ and $O_{\delta}$ be as in the proof of Theorem 1.1. Given
  $\rho\in O_{\delta}$, we first solve the equation $(2.11)_1$ subject to the boundary value condition $(2.11)_3$.
  By Lemma 2.4 and Theorem 1.1, this problem has a unique solution $\sigma$. Next we substitute $\sigma$ into $(2.11)_2$
  and take $(2.11)_4$ into account. Then we obtain the following elliptic boundary value problem:
\begin{equation}
\left\{
\begin{array}{rll}
   -\Delta\pi_0=&g(\sigma) &\quad\;\; \mbox{in} \;\; \Omega,\\
   \pi_0=&0   &\quad\;\; \mbox{on} \;\; \partial\Omega.
\end{array}
\right.
\end{equation}
  Since $g(\sigma)\in L^{\infty}(\Omega)$, by applying standard theory for elliptic boundary value problems and using the
  properties of $\sigma$ proved in Lemma 2.4, we see the above problem has a unique solution satisfying the following properties:
\begin{equation}
  \pi_0\in W^{2,q}(\Omega) \;(\forall q\in [1,\infty)), \quad
  \pi_0\in C^{\infty}(\Omega_{\rm liv}\cup\Omega_{\rm nec}) \quad \mbox{and} \quad
  \pi_0\in C^{m+\mu}(\Omega_{\rm liv}\cup\partial\Omega).
\end{equation}
  It follows that $\partial_{\bfn}\pi|_{\partial\Omega}\in C^{m-1+\mu}(\partial\Omega)$, where $\bfn$ denotes the unit
  outward normal field of $\partial\Omega$. In this way we obtain a map $F_0:O_{\delta}\subseteq C^{m+\mu}({\mathbf{S}}^2)
  \to C^{m-1+\mu}({\mathbf{S}}^2)$ defined as follows: For any $\rho\in O_{\delta}$,
$$
  F_0(\rho)=[\omega\mapsto\partial_{\bfn}\pi_0(R[1+\rho(\omega)],\omega),\omega\in{\mathbf{S}}^2].
$$
  Our next goal of this section is to prove the following result:
\medskip

  {\bf Lemma 2.5}\ \ {\em $F_0\in C^{\infty}(O_{\delta},C^{m-1+\mu}({\mathbf{S}}^2))$.}
\medskip

  {\em Proof}.\ \ From the proof of Lemma 2.2 we see that not only the map $\rho\mapsto\eta$ is tame, but also
  the map $\rho\mapsto\sigma|_{\overline{E}}$ from  $O_{\delta}\subseteq C^{m+\mu}({\mathbf{S}}^2)$ to
  $C^{\infty}(\overline{E})$ is tame. Let $D_{\rho,\eta}$, $S_{\rho}$ and $\Gamma_{\eta}$ be as before and set
$$
  B_{\eta}=\{x\in{\mathbf{R}}^3: r< K[1+\eta(\omega)]\}, \qquad B_0=B(0,K).
$$
  Then (2.13) can be rewritten as the following equivalent problem:
\begin{equation}
\left\{
\begin{array}{rll}
   -\Delta\pi_0=&a(\sigma-\tilde{\sigma})-b & \;  \mbox{in}\;\; D_{\rho,\eta},\\
   -\Delta\pi_0=&-b & \;  \mbox{in}\;\; B_{\eta},\\
   \pi_0=&0 & \;  \mbox{on}\;\; S_{\rho},\\
   \pi_0,\partial_{\nu}\pi_0\;& \;\mbox{are continuous}\;& \; \mbox{across}\;\, \Gamma_{\eta},
\end{array}
\right.
\end{equation}
  where $\nu$ is as before (note that it is also the inward unit normal field of the boundary $\Gamma_{\eta}$
  of $B_{\eta}$). Let $\Psi_{\rho,\eta}$, $\Psi_{\eta}$, $\psi_{\eta}$, $\mathscr{A}(\rho,\eta)$, $u$ be as in
  the proof of Lemma 2.2 and define
$$
  \mathscr{B}(\eta)u=[\partial_{\nu}(u\circ\Psi_{\eta})|_{\Gamma_{\eta}}]\circ\psi_{\eta}^{-1} \quad\;\;
  \mbox{for}\;\,u\in C^{\infty}(\overline{E}).
$$
  Choose another smooth function $\phi_1\in C^{\infty}[0, K]$ such that it satisfies the following conditions:
$$
  0\leqslant\phi_1\leqslant 1; \quad \phi_1'\geqslant 0; \quad \phi_1(t)=0 \;\; \mbox{for}\;\,0\leqslant t
  \leqslant\frac{1}{2} K; \quad
  \phi_1( K)=1.
$$
  Let $M_1=\displaystyle\max_{0\leqslant t\leqslant K}|\phi_1'(t)|$ and assume $\delta'>0$ is small enough such that
  in addition to the conditions appearing in the proof of Lemma 2.2 we have also $\delta'<(1+M_1 K)^{-1}$. Let
  $\Psi_{\eta}^1:\overline{B}_{\eta}\to \overline{B}_0$ be as follows:
\begin{equation}
  \Psi_{\eta}^1(x)=x- K\eta(\omega)\phi_1\Big(\frac{r}{1+\eta(\omega)}\Big)\omega \quad\;\;
  \mbox{for} \;\; x\in\overline{B}_{\eta}.
\end{equation}
  Define $\mathscr{A}_1(\eta):C^{\infty}(\overline{B}_0)\to C^{\infty}(\overline{B}_0)$ and $\mathscr{B}_1(\eta):
  C^{\infty}(\overline{B}_0)\to C^{\infty}(\Gamma_0)$ respectively as follows:
$$
  \mathscr{A}_1(\eta)u=[\Delta(u\circ\Psi_{\eta}^1)]\circ(\Psi_{\eta}^1)^{-1} \quad\;\; \mbox{for}\;\,u\in
  C^{\infty}(\overline{B}_0),
$$
$$
  \mathscr{B}_1(\eta)u=[\partial_{\nu}(u\circ\Psi_{\eta}^1)|_{\Gamma_{\eta}}]\circ(\psi_{\eta}^1)^{-1} \quad\;\;
  \mbox{for}\;\,u\in C^{\infty}(\overline{B}_0),
$$
  where $\psi_{\eta}^1=\Psi_{\eta}^1|_{\Gamma_{\eta}}$. Let $v=\pi_0|_{\overline{D}_{\rho,\eta}}\circ\Psi_{\rho,\eta}^{-1}$
  and $v_1=\pi_0|_{\overline{B}_{\eta}}\circ(\Psi_{\eta}^1)^{-1}$. After the variable transformation $x\mapsto
  \Psi_{\rho,\eta}(x)$ (for $x\in\overline{D}_{\rho,\eta}$) and $x\mapsto\Psi_{\eta}^1(x)$ (for $x\in B_{\eta}$),
  the problem (2.15) transforms into the following problem:
\begin{equation}
\left\{
\begin{array}{rll}
   -\mathscr{A}(\rho,\eta)v=&a(u-\tilde{\sigma})-b &  \quad  \mbox{in}\;\; D,\\
   -\mathscr{A}_1(\eta)v_1=&-b &  \quad  \mbox{in}\;\; B_0,\\
   v=&0 &  \quad  \mbox{on}\;\; S_0,\\
   v=&v_1 &  \quad  \mbox{on}\;\; \Gamma_0,\\
   \mathscr{B}(\eta)v=& \mathscr{B}_1(\eta)v_1 &   \quad  \mbox{on}\;\; \Gamma_0.
\end{array}
\right.
\end{equation}
  Lemma 2.5 is an immediate consequence of Theorem 1.1 and the following preliminary result:
\medskip

  {\bf Lemma 2.6}\ \ {\em Let $\delta,O_{\delta}$ be as in the proof of Lemma 2.2 and $\delta'$ as above.
  Let $O_{\delta'}''=\{\eta\in C^{m+\mu}({\mathbf{S}}^2): \|\eta\|_{C^{m+\mu}({\mathbf{S}}^2)}<\delta'\}$. Given
  $(\rho,\eta,u)\in O_{\delta}\times O_{\delta'}''\times C^{m+\mu}(\overline{D})$, the problem $(2.17)$ has a unique
  solution $(v,v_1)\in C^{m+\mu}(\overline{D})\times C^{m+\mu}(\overline{B}_0)$, and the solution map $(\rho,\eta,u)
  \mapsto (v,v_1)$ from $O_{\delta}\times O_{\delta'}''\times C^{m+\mu}(\overline{D})\subseteq C^{m+\mu}({\mathbf{S}}^2)
  \times C^{m+\mu}({\mathbf{S}}^2)\times C^{m+\mu}(\overline{D})$ to $C^{m+\mu}(\overline{D})\times
  C^{m+\mu}(\overline{B}_0)$ is smooth.}
\medskip

  {\em Proof.}\ \ We first note that for $\Omega=B(0,R)$ with $R>R^*$, the unique solution of the problem (2.13) is
  given by $\pi_0=V(r,R)$, where
$$
  V(r,R)=\left\{
\begin{array}{ll}
    \displaystyle D\Big(\frac{1}{R}\!-\!\frac{1}{r}\Big)-\frac{1}{6}(a\tilde{\sigma}+b)(R^2\!-\!r^2)
    -a\int_r^R\!\!\!\int_{\xi}^RU(\eta,R)\Big(\frac{\eta}{\xi}\Big)^2d\eta d\xi &\quad\;\;
    \mbox{for} \;\;  K\leqslant r\leqslant R,\\ [0.3cm]
   \displaystyle C+\frac{b}{6}r^2  &\quad\;\; \mbox{for} \;\; r< K,
\end{array}
\right.
$$
  where $U(r,R)$ and $K=K(R)$ are as before, and $C,D$ are constants such that the relations $V( K^+,R)=V( K^-,R)$
  and $\partial_rV( K^+,R)=\partial_rV( K^-,R)$ hold; in particular,
$$
  D=-\frac{1}{3}a\tilde{\sigma}K^3-a\int_{ K}^RU(\eta,R)\eta^2d\eta.
$$
  (cf. Lemmas 3.1 and 3.2 of \cite{Cui2}).

  Given $(\rho,\eta,u,\xi)\in O_{\delta}\times O_{\delta'}''\times C^{m+\mu}(\overline{D})\times C^{m+\mu}(\Gamma_0)$, we
  consider the following two elliptic boundary value problems:
\begin{equation}
\left\{
\begin{array}{rll}
   -\mathscr{A}(\rho,\eta)v=\,&a(u-\tilde{\sigma})-b &  \quad  \mbox{in}\;\; D,\\
   v=\,&0 &  \quad  \mbox{on}\;\; S_0,\\
   v=\,&\xi &  \quad  \mbox{on}\;\; \Gamma_0,
\end{array}
\right.
  \qquad\qquad
\left\{
\begin{array}{rll}
   -\mathscr{A}_1(\eta)v_1=\,&-b &  \quad  \mbox{in}\;\; B_0,\\
   v_1=\,&\xi &  \quad  \mbox{on}\;\; \Gamma_0.
\end{array}
\right.
\end{equation}
  Clearly, these problems have unique solutions $v\in C^{m+\mu}(\overline{D})$ and $v_1\in C^{m+\mu}(\overline{B}_0)$,
  respectively. Define $\mathscr{D}(\rho,\eta):C^{m+\mu}(\Gamma_0)\to C^{m-1+\mu}(\Gamma_0)$ and $\mathscr{D}_1(\eta):
  C^{m+\mu}(\Gamma_0)\to C^{m-1+\mu}(\Gamma_0)$ respectively as follows:
$$
  \mathscr{D}(\rho,\eta)\xi=\mathscr{B}(\eta)v, \quad \mathscr{D}_1(\eta)\xi=\mathscr{B}_1(\eta)v_1 \quad
  \mbox{for}\;\,\xi\in C^{m+\mu}(\Gamma_0).
$$
  The problem (2.17) is equivalent to the following problem: Find $\xi\in C^{m+\mu}(\Gamma_0)$ such that
\begin{equation}
  \mathscr{D}(\rho,\eta)\xi=\mathscr{D}_1(\eta)\xi.
\end{equation}
  We introduce a mapping $\mathscr{G}: O_{\delta}\times O_{\delta'}''\times C^{m+\mu}(\overline{D})\times
  C^{m+\mu}(\Gamma_0)\to C^{m-1+\mu}(\Gamma_0)$ by defining
$$
  \mathscr{G}(\rho,\eta,u,\xi)=\mathscr{D}(\rho,\eta)\xi-\mathscr{D}_1(\eta)\xi \quad
  \mbox{for}\;\,\rho\in O_{\delta},\;\,\eta\in O_{\delta'}'',\;\, u\in C^{m+\mu}(\overline{D}),\;\,
  \xi\in C^{m+\mu}(\Gamma_0).
$$
  It is clear that $\mathscr{G}\in C^{\infty}(O_{\delta}\times O_{\delta'}'' C^{m+\mu}(\overline{D})\times
  C^{m+\mu}(\Gamma_0),C^{m-1+\mu}(\Gamma_0))$, and
$$
  \mathscr{G}(0,0,u_0,\xi_0)=0,
$$
  where $u_0=U(r,R)$ and $\xi_0$ represents the following constant function in $\Gamma_0$: $\xi_0(x)=V( K,R)$ for
  $x\in\Gamma_0$. A simple computation shows that for any $\zeta\in C^{m+\mu}(\Gamma_0)$, $L\zeta:=\partial_{\xi}
  \mathscr{G}(0,0,u_0,\xi_0)\zeta=\partial_{\nu}w|_{\Gamma_0}-\partial_{\nu}w_1|_{\Gamma_0}=\partial_rw_1|_{\Gamma_0}
  -\partial_rw|_{\Gamma_0}$, where $w,w_1$ are respectively the unique solutions of the following problems:
\begin{equation}
\left\{
\begin{array}{rll}
   \Delta w=&0 &  \quad  \mbox{for}\;\;  K<r<R,\\
   w=&0 &  \quad  \mbox{for}\;\; r=R,\\
   w=&\zeta &  \quad  \mbox{for}\;\; r=K,
\end{array}
\right.
  \qquad\qquad
\left\{
\begin{array}{rll}
   \Delta w_1=&0 &  \quad  \mbox{for}\;\; r< K,\\
   w_1=&\zeta &  \quad  \mbox{for}\;\; r= K.
\end{array}
\right.
\end{equation}
  From this fact it is not hard to see that if $L\zeta=0$ then $\zeta=0$. Indeed, if $L\zeta=0$ then by letting $z=w$
  for $ K\leqslant r\leqslant R$ and $z=w_1$ for $r< K$, we get a weak solution of the boundary value problem
$$
\left\{
\begin{array}{rl}
   \Delta z=&0   \quad  \mbox{in}\;\; B(0,R),\\
   z=&0   \quad  \mbox{on}\;\;\partial B(0,R),
\end{array}
\right.
$$
  which implies, by Green's second identity, that $z=0$ and, consequently, $\zeta=0$. This shows that
  ${\rm Ker}\,L=\{0\}$. Since $L$ is a sum of two Dirichlet-Neumann operators, it is a first-order pseudo-differential
  operator of elliptic type (cf. \cite{EscSim2},  \cite{EscSim5}), so that standard Schauder estimate for elliptic
  pseudo-differential operator applies to it: There exists a positive constant $C>0$ such that
$$
  \|\zeta\|_{C^{m+\mu}(\Gamma_0)}\leqslant C(\|L\zeta\|_{C^{m-1+\mu}(\Gamma_0)}+\|\zeta\|_{L^{\infty}(\Gamma_0)}),
   \quad \forall\zeta\in C^{m+\mu}(\Gamma_0).
$$
   Since ${\rm Ker}\,L=\{0\}$, this implies, by a standard argument, the following estimate:
\begin{equation}
  \|\zeta\|_{C^{m+\mu}(\Gamma_0)}\leqslant C\|L\zeta\|_{C^{m-1+\mu}(\Gamma_0)}, \quad \forall\zeta\in C^{m+\mu}(\Gamma_0).
\end{equation}
  For every $k\in{\mathbf{Z}}_+$ let $\{Y_{kl}(\omega)\}_{l=1}^{2k+1}$ be the normalized orthogonal basis
  (in $L^2({\mathbf{S}}^2)$ inner product) of the linear space of $k$-th order spherical harmonics. A simple computation
  shows that for any $\zeta\in C^{\infty}(\Gamma_0)$,
$$
  L\zeta(\omega)=\sum_{k=0}^{\infty}\sum_{l=1}^{2k+1}\frac{(2k\!+\!1)c_{kl}}{[1-(K/R)^{2k+1}]K}Y_{kl}(\omega) \quad
  \mbox{if}\;\; \zeta(K\omega)=\sum_{k=0}^{\infty}\sum_{l=1}^{2k+1}c_{kl}Y_{kl}(\omega).
$$
  From this expression of $L$ it is easy to see that for any $\eta\in C^{\infty}(\Gamma_0)$ the equation $L\zeta=\eta$
  has a unique solution $\zeta\in C^{\infty}(\Gamma_0)$ (see the proofs of Lemma 3.3 and Corollary 3.4 in the next section
  for more details in this argument). It follows, by using the estimate (2.21) and a standard approximation argument,
  that also for any $\eta\in C^{m-1+\mu}(\Gamma_0)$ the equation $L\zeta=\eta$ has a unique solution $\zeta\in
  C^{m+\mu}(\Gamma_0)$\footnotemark[4].
\footnotetext[4]{Since $C^{\infty}(\Gamma_0)$ is not dense in $C^{m-1+\mu}(\Gamma_0)$, one might argue validity of the
  approximation argument. It is as follows: For any $\eta\in C^{m-1+\mu}(\Gamma_0)$ choose a sequence
  $\{\eta_j\}_{j=1}^{\infty}\subseteq C^{\infty}(\Gamma_0)$ such that it is bounded in $C^{m-1+\mu}(\Gamma_0)$ and converges
  to $\eta$ in $C^{m-1+\mu'}(\Gamma_0)$ for any $0<\mu'<\mu$ (e.g., choose a mollification sequence of $\eta$). For every
  $j$ let $\zeta_j\in C^{\infty}(\Gamma_0)$ be the unique solution of the equation $L\zeta_j=\eta_j$. The estimate (2.21)
  ensures the sequence $\{\zeta_j\}_{j=1}^{\infty}$ is bounded in $C^{m+\mu}(\Gamma_0)$ and converges to a function $\zeta$
  in $C^{m+\mu'}(\Gamma_0)$ for any $0<\mu'<\mu$, which implies that $\zeta$ is a solution of the equation $L\zeta=\eta$.
  Boundedness of the sequence $\{\zeta_j\}_{j=1}^{\infty}$ in $C^{m+\mu}(\Gamma_0)$ implies $\zeta\in C^{m+\mu}(\Gamma_0)$.}
  This shows that $L=\partial_{\xi}\mathscr{G}(0,0,u_0,\xi_0): C^{m+\mu}(\Gamma_0)\to C^{m-1+\mu}(\Gamma_0)$ is a (linear
  and topological) isomorphism. Hence, by applying the implicit function theorem in Banach spaces, we see that by choosing
  $\delta,\delta'>0$ further small when necessary, there exists a smooth mapping $\chi:O_{\delta}\times O_{\delta'}''\times
  C^{m+\mu}(\overline{D})\subseteq C^{m+\mu}({\mathbf{S}}^2)\times C^{m+\mu}({\mathbf{S}}^2)\times C^{m+\mu}(\overline{D})
  \to C^{m+\mu}(\Gamma)$ such that $\chi(0,0,u_0)=\xi_0$ and for any $\rho\in O_{\delta}$, $\eta\in O_{\delta'}''$ and
  $u\in C^{m+\mu}(\overline{D})$, $\xi=\chi(\rho,\eta,u)$ is the unique solution of the equation $\mathscr{G}(\rho,\eta,u,\xi)
  =0$ in a small neighborhood of $\xi_0$ in $C^{m+\mu}(\Gamma)$. This proves unique solvability of the equation (2.19)
  and smoothness of the solution map $\xi=\chi(\rho,\eta,u)$. As a result, the solution map $(\rho,\eta,u)\mapsto(v,v_1)$
  (from $O_{\delta}\times O_{\delta'}''\times C^{m+\mu}(\overline{D})\subseteq C^{m+\mu}({\mathbf{S}}^2)\times
  C^{m+\mu}({\mathbf{S}}^2)\times C^{m+\mu}(\overline{D})$ to $C^{m+\mu}(\overline{D})\times C^{m+\mu}(\overline{B}))$) of
  the problem (2.17) is smooth. This proves Lemma 2.6 and also completes the proof of Lemma 2.5. $\quad\Box$
\medskip

   {\em Remark}.\ \ We note that all H\"{o}lder spaces appearing in this section can be replaced with corresponding
   little H\"{o}lder spaces. In the  next section we shall use this fact without making further explanation.

\section{The proof of Theorem 1.2}
\setcounter{equation}{0}

\hskip 2em
  In this section we prove Theorem 1.2.

  We first recall that radial stationary solution $(\sigma_s,\pi_s,\Omega_s)$ of the problem (1.1) is given by
\begin{equation}
  \sigma_s(r)=U(r,R_s), \qquad  \pi_s(r)=\frac{\gamma}{R_s}+V(r,R_s), \qquad \Omega_s=B(0,R_s),
\end{equation}
  where $U$, $V$ are as in the previous section, and $R_s$ is the root of the equation $\pi_s'(R_s)=0$ or
  $\partial_rV(R_s,R_s)=0$, i.e.,
$$
  \frac{1}{3}a\tilde{\sigma}K^3(R_s)+a\int_{K(R_s)}^{R_s}U(\eta,R_s)\eta^2d\eta=\frac{1}{3}(a\tilde{\sigma}+b)R_s^3.
$$
  Since $0<\tilde{\sigma}<1$, by Lemma 4.2 of \cite{Cui2} we know that this equation has a unique solution $R_s>R^*$,
  which means the problem (1.1) with $f$, $g$ given in (1.2) has a unique radial stationary solution.

  We shall use Theorem 1.1 of \cite{Cui5} to prove Theorem 1.2. To this end, let us first reduce the problem (1.1)
  into a differential equation in a Banach manifold. Let $m$ be a positive integer$\,\geqslant$2 and $0<\mu<1$. Let
  $\mathfrak{M}:=\dot{\mathfrak{D}}^{m+\mu}({\mathbf{R}}^3)$ and $\mathfrak{M}_0:=
  \dot{\mathfrak{D}}^{m+\mu+3}({\mathbf{R}}^3)$. Given $\Omega\in\mathfrak{M}_0$, we have seen that the problem (2.11)
  has a unique solution $(\sigma,\pi_0)$ satisfying the following properties:
$$
  \sigma,\pi_0\in W^{2,q}(\Omega) \;(\forall q\in [1,\infty)), \quad
  \sigma,\pi_0\in\dot{C}^{\infty}(\Omega_{\rm liv}\cup\Omega_{\rm nec}) \quad \mbox{and} \quad
  \sigma,\pi_0\in\dot{C}^{m+3+\mu}(\Omega_{\rm liv}\cup\partial\Omega).
$$
  It follows that $\partial_{\bfn}\pi_0|_{\partial\Omega}\in \dot{C}^{m+\mu+2}(\partial\Omega)$.
  We define $\mathscr{F}_0:\mathfrak{M}_0\to\mathcal{T}_{\mathfrak{M}_0}(\mathfrak{M})$ by letting
$$
  \mathscr{F}_0(\Omega)=-\partial_{\bfn}\pi_0|_{\partial\Omega}, \quad  \forall\Omega\in\mathfrak{M}_0.
$$
  Next, given $\Omega\in\mathfrak{M}_0$, let $\pi\in\dot{C}^{m+1+\mu}(\overline{\Omega})$ be the unique solution of
  the following elliptic boundary value problem:
$$
\left\{
\begin{array}{rll}
   -\Delta\pi=&0   &\quad\;\; \mbox{in} \;\; \Omega,\\
   \pi=&\kappa  &\quad\;\; \mbox{on} \;\; \partial\Omega,
\end{array}
\right.
$$
  where $\kappa$ is as explained in Section 1, and define $\mathscr{F}:\mathfrak{M}_0\to
  \mathcal{T}_{\mathfrak{M}_0}(\mathfrak{M})$ by letting
$$
  \mathscr{F}(\Omega)=-\partial_{n}\pi|_{\partial\Omega}, \quad  \forall\Omega\in\mathfrak{M}_0.
$$
  We now define $\mathscr{G}:\mathfrak{M}_0\to\mathcal{T}_{\mathfrak{M}_0}(\mathfrak{M})$ as follows:
\begin{equation}
  \mathscr{G}(\Omega)=\gamma\mathscr{F}(\Omega)+\mathscr{F}_0(\Omega), \quad  \forall\Omega\in\mathfrak{M}_0.
\end{equation}
  Then $\mathscr{F}$ is a vector field in $\mathfrak{M}$ with domain $\mathfrak{M}_0$, and the problem (1.1) reduces
  into the following differential equation in the Banach manifold $\mathfrak{M}$:
\begin{equation}
\left\{
\begin{array}{ll}
   \Omega'(t)=\mathscr{G}(\Omega(t)), &\quad  t>0,\\
    \Omega(0)=\Omega_0. &
\end{array}
\right.
\end{equation}
  The fact that $(\sigma_s,\pi_s,\Omega_s)$ is a stationary solution of the problem (1.1) implies that $\Omega_s$ is
  a stationary solution of the equation $\Omega'=\mathscr{G}(\Omega)$.

  Let $G_{tl}={\mathbb{R}}^n$ be the additive group of $n$-vectors. Given $z\in{\mathbb{R}}^n$ and $\Omega\in
  \mathfrak{M}$, let
$$
  p(z,\Omega)=\Omega+z=\{x+z:\,x\in \Omega\}.
$$
  It is clear that $p(z,\Omega)\in\mathfrak{M}$, $\forall\Omega\in\mathfrak{M}$, $\forall z\in{\mathbb{R}}^n$. It can
  be easily seen that $(G_{tl},p)$ is a Lie group action on $\mathfrak{M}$. By Lemma 4.1 of \cite{Cui3} we know that
  the action $p(z,\Omega)$ is differentiable at every point $\Omega\in\mathfrak{M}_0$, and $\mbox{rank}D_{z}p(z,\Omega)=n$,
  $\forall z\in G_{tl}$, $\forall S\in\mathfrak{M}_0$.
\medskip

  {\bf Lemma 3.1}\ \ {\em The vector field $\mathscr{G}$ is invariant under the group action $(G_{tl},p)$.}
\medskip

  {\em Proof}.\ \ The proof is similar to that of Lemma 5.3 of \cite{Cui5}. We omit it here. $\quad\Box$
\medskip

  Next we consider representation of the problem (3.3) in a regular local chart $(\mathcal{U},\varphi)$
  of $\mathfrak{M}$ at the point $\Omega_s$. We denote
$$
  X=\dot{C}^{m+\mu}({\mathbf{S}}^2),  \qquad  X_0=\dot{C}^{m+\mu+3}({\mathbf{S}}^2),  \qquad  X_1=\dot{C}^{m+\mu+2}({\mathbf{S}}^2),
$$
  and for a sufficiently small number $\delta>0$ let
$$
  O_{\delta}=\{\rho\in X,\|\rho\|_{X}<\delta\},  \qquad O'_{\delta}=\{\rho\in X_0,\|\rho\|_{X_0}<\delta\},
$$
$$
  \mathcal{U}=\{\Omega_{\rho}:\rho\in O_{\delta}\},  \qquad \mathcal{U}'=\{\Omega_{\rho}:\rho\in O'_{\delta}\},
$$
  where $\Omega_{\rho}$ is as before, i.e.,
$$
  \Omega_{\rho}=\{x\in{\mathbf{R}}^3: r<R_s[1+\rho(\omega)]\}.
$$
  It is clear that $\Omega_{\rho}\in\mathfrak{M}$ for $\rho\in O_{\delta}$ and $\Omega_{\rho}\in\mathfrak{M}_0$ for
  $\rho\in O_{\delta}'$, so that $\mathcal{U}$ and $\mathcal{U}'$ are neighborhoods of $\Omega_s$ in $\mathfrak{M}$
  and $\mathfrak{M}_0$, respectively. We define $\varphi:\mathcal{U}\to X$ by letting
  $\varphi(\Omega_{\rho})=\rho$, $\forall\rho\in O_{\delta}$. Then $(\mathcal{U},\varphi)$ is a regular local chart
  of $\mathfrak{M}$ at the point $\Omega_s$, with base space $X$. We denote by $F$, $F_0$ and $G$ the representations
  of the vector fields $\mathscr{F}$, $\mathscr{F}_0$ and $\mathscr{G}$, respectively, in this local chart, i.e.,
  for any $\rho\in O'_{\delta}$,
$$
  F(\rho)=\varphi'(\Omega_{\rho})\mathscr{F}(\Omega_{\rho}), \quad F_0(\rho)=\varphi'(\Omega_{\rho})
  \mathscr{F}_0(\Omega_{\rho})
   \quad \mbox{and} \quad G(\rho)=\varphi'(\Omega_{\rho})\mathscr{G}(\Omega_{\rho});
$$
  see Section 2 of \cite{Cui5} for details of these notions.
  Then $G(\rho)=\gamma F(\rho)+F_0(\rho)$, and representation of the problem (3.3) in the local chart
  $(\mathcal{U},\varphi)$ is the following initial value problem in the Banach space $X$:
\begin{equation}
\left\{
\begin{array}{ll}
   \rho'(t)=G(\rho(t)), &\quad  t>0,\\
   \rho(0)=\rho_0, &
\end{array}
\right.
\end{equation}
  where $\rho_0\in O'_{\delta}$ is the function such that $\Omega_0=\Omega_{\rho_0}$. It is clear that $F\in
  C^{\infty}(O'_{\delta},X)$. By Lemma 2.5 we know that also $F_0\in C^{\infty}(O'_{\delta},X_1)\subseteq
  C^{\infty}(O'_{\delta},X)$. Hence we have
\medskip

  {\bf Lemma 3.2}\ \ {\em $G\in C^{\infty}(O'_{\delta},X)$.} $\quad\Box$
\medskip

  {\bf Lemma 3.3}\ \ {\em The differential equation $(3.3)$ is of parabolic type.}
\medskip

  {\em Proof.}\ \ It is well-known that for any $\rho\in O'_{\delta}$, $F'(\rho)$ is a sectorial operator in $X$
  with domain $X_0$; cf., e.g., \cite{Esc}, \cite{Gri}. Lemma 2.5 ensures that for any $\rho\in O'_{\delta}$,
  $F_0'(\rho)\in L(X_0,X_1)$. Since
  $G'(\rho)=\gamma F'(\rho)+F_0'(\rho)$ and $X_1$ is an intermediate space between $X_0$ and $X$, by a well-known
  perturbation theorem for sectorial operators it follows that $G'(\rho)$ is also a sectorial operator in $X$ with
  domain $X_0$. Besides, it is easy to check that the graph norm of ${\rm Dom}\,G'(\rho)=X_0$ is equivalent to the
  norm of $X_0$. Hence the desired assertion follows. $\quad\Box$

  Following \cite{FriRei2} and \cite{CuiEsc2}, we compute $G'(0)$ as follows: Let
$$
  \rho(\omega)=\varepsilon\xi(\omega), \quad \eta(\omega)=\varepsilon\zeta(\omega), \quad \sigma(r,\omega)=\sigma_s(r)+\varepsilon u(r,\omega).
$$
  Substituting these expressions into (2.12) and the equation $\Delta\sigma=0$ for $r< K_s$, and comparing coefficients of first-order terms of
  $\varepsilon$, we obtain the following equations:
\begin{equation}
\left\{
\begin{array}{rll}
   \Delta u=&u &  \quad  \mbox{for}\;\;  K_s<r<R_s,\\
   \Delta u=&0 &  \quad  \mbox{for}\;\; r< K_s,\\
   u=&-R_s\sigma_s'(R_s)\xi(\omega) &  \quad  \mbox{for}\;\; r=R_s,\\
   u=&0 &  \quad  \mbox{for}\;\; r= K_s,\\
   \partial_r^+u=&-\hat{\sigma} K_s\zeta(\omega) &  \quad  \mbox{for}\;\; r= K_s.
\end{array}
\right.
\end{equation}
  Here $K_s=K(R_s)$ and for $r= K_s$, $\partial_r^+u=\frac{\partial u}{\partial r}( K_s^+,\omega)$. Similarly,
  by letting $\pi(r,\omega)=\pi_s(r)+\varepsilon v(r,\omega)$, we get the following equations for $v$:
\begin{equation}
\left\{
\begin{array}{rll}
   \Delta v=&-au &  \quad  \mbox{for}\;\;  K_s<r<R_s,\\
   \Delta v=&0 &  \quad  \mbox{for}\;\; r< K_s,\\
   v=&\displaystyle-\frac{\gamma}{R_s}\Big(\xi(\omega)+\frac{1}{2}\Delta_{\omega}\xi(\omega)\Big) &
   \quad  \mbox{for}\;\; r=R_s,\\
   v^+=&v^- &  \quad  \mbox{for}\;\; r= K_s,\\
   \partial_r^+v-\partial_r^-v=&a(\hat{\sigma}-\tilde{\sigma}) K_s\zeta(\omega) &  \quad  \mbox{for}\;\; r= K_s.
\end{array}
\right.
\end{equation}
  Here for $r=K_s$, $v^{\pm}=v( K_s^{\pm},\omega)$ and $\partial_r^{\pm}v=\frac{\partial v}{\partial r}( K_s^{\pm},\omega)$.
  After solving these equations (with $u$, $v$ and $\zeta$ being unknown functions) for any given function $\xi=\xi(\omega)$,
  we then have
\begin{equation}
  G'(0)\xi(\omega)=-\frac{\partial v}{\partial r}(R_s,\omega)+g(1)R_s\xi(\omega), \quad \forall\xi\in X_0
\end{equation}
  (cf. (4.4) in \cite{CuiEsc2}, but be aware that here the perturbation of the sphere $r=R_s$ is given by
  $r=R_s[1+\varepsilon\xi(\omega)]$, not as in \cite{CuiEsc2} given by $r=R_s+\varepsilon\xi(\omega)$).

  We now use spherical harmonics expansions of functions in ${\mathbf{S}}^2$ to solve the problems (3.5) and (3.6).
  Hence let $\{Y_{kl}(\omega): k=0,1,2,\cdots,\; l=1,2,\cdots,2k+1\}$ be a normalized orthogonal basis of
  $L^2({\mathbf{S}}^2)$ consisting spherical harmonics on ${\mathbf{S}}^2$, where for every $k\in{\mathbf{Z}}_+$,
  $Y_{kl}(\omega)$, $l=1,2,\cdots,2k+1$, are spherical harmonics of degree $k$, so that
$$
  \Delta_{\omega}Y_{kl}(\omega)=-\lambda_kY_{kl}(\omega), \quad k=0,1,2,\cdots,\;\; l=1,2,\cdots,2k+1.
$$
  where $\lambda_k=k(k+1)$, $k=0,1,2,\cdots$. A simple computation shows that if a given function
  $\xi\in C^{\infty}({\mathbf{S}}^2)$ has a spherical harmonics expansion
\begin{equation}
  \xi(\omega)=\sum_{k=0}^{\infty}\sum_{l=1}^{2k+1}c_{kl}Y_{kl}(\omega),
\end{equation}
  then the solution $(u,\zeta)$ of the problem (3.5) is given by
$$
  u(r,\omega)=\left\{
\begin{array}{cl}
  \displaystyle -R_s\sigma_s(R_s)\sum_{k=0}^{\infty}\sum_{l=1}^{2k+1}\Big(\frac{r}{R_s}\Big)^k\bar{u}_k(r)c_{kl}Y_{kl}(\omega) &  \quad
   \mbox{for}\;\;  K_s<r\leqslant R_s\\
   0 &  \quad  \mbox{for}\;\; r\leqslant K_s
\end{array}
\right.
$$
  and $\zeta(\omega)=-\partial_r^+u/\hat{\sigma} K_s$, where $\bar{u}_k$ is the unique solution of the following boundary value problem:
\begin{equation}
  \left\{
\begin{array}{l}
  \displaystyle \bar{u}_k''(r)+\frac{2(k\!+\!1)}{r}\bar{u}_k'(r)=\bar{u}_k(r) \quad  \mbox{for}\;\;  K_s<r<R_s,\\
   \bar{u}_k( K_s)=0,   \quad  \bar{u}_k(R_s)=1.
\end{array}
\right.
\end{equation}
  Substituting these expressions of $u$ and $\zeta$ into (3.6), we easily obtain the following solution of that problem:
\begin{equation}
  v(r,\omega)=\sum_{k=0}^{\infty}\sum_{l=1}^{2k+1}\Big[\frac{\gamma}{2R_s}(k-1)(k+2)
  +aR_s\sigma_s'(R_s)\bar{v}_k(r)\Big]\Big(\frac{r}{R_s}\Big)^kc_{kl}Y_{kl}(\omega),
\end{equation}
  where $\bar{v}_k$ is the unique solution of the following boundary value problem:
\begin{equation}
  \left\{
\begin{array}{l}
  \displaystyle \bar{v}_k''(r)+\frac{2(k\!+\!1)}{r}\bar{v}_k'(r)=\bar{u}_k(r) \quad  \mbox{for}\;\;  K_s<r<R_s,\\
  \displaystyle \bar{v}_k''(r)+\frac{2(k\!+\!1)}{r}\bar{v}_k'(r)=0  \quad \mbox{for}\;\; r< K_s,\\
   \bar{v}_k(R_s)=0,   \\
   \bar{v}_k( K_s^+)=\bar{v}_k( K_s^-),  \\
  \displaystyle \bar{v}_k'( K_s^+)-\bar{v}_k'( K_s^-)=\frac{\hat{\sigma}-\tilde{\sigma}}{\hat{\sigma}}\bar{u}_k'( K_s^+).
\end{array}
\right.
\end{equation}
  From (3.11) we get the following relation:
\begin{equation}
  \bar{v}_k'(R_s)=\frac{\hat{\sigma}-\tilde{\sigma}}{\hat{\sigma}}\bar{u}_k'(K^+)\Big(\frac{K_s}{R_s}\Big)^{2(k+1)}
  +\int_{K_s}^{R_s}\bar{u}_k(\tau)\Big(\frac{\tau}{R_s}\Big)^{2(k+1)}d\tau.
\end{equation}
  Now let
\begin{equation}
  a_k(\gamma)=-\frac{\gamma}{2R_s^2}k(k-1)(k+2)-aR_s\sigma_s'(R_s)\bar{v}_k'(R_s)+g(1)R_s.
\end{equation}
  Then from (3.8), (3.10), $(3.11)_3$ and (3.12) we obtain the following preliminary result:
\medskip

  {\bf Lemma 3.4}\ \ {\em $G'(0)$ is a Fourier multiplier in the sense that if $\xi\in C^{\infty}({\mathbf{S}}^2)$ has
  expansion $(3.8)$, then $G'(0)\xi=\displaystyle\sum_{k=0}^{\infty}\sum_{l=1}^{2k+1}a_k(\gamma)c_{kl}Y_{kl}(\omega)\in
  C^{\infty}({\mathbf{S}}^2)$.}  $\quad\Box$
\medskip

  {\bf Corollary 3.5}\ \ {\em $\sigma(G'(0))=\{a_k(\gamma):k=0,1,2,\cdots\}$.}
\medskip

  {\em Proof.}\ \ Since $G'(0)\in L(X_0,X)$ and it is a sectorial operator in $X$ with domain $X_0$, the inverse mapping theorem implies that for
  any $\lambda\in\rho(G'(0))$ we have $[\lambda I-G'(0)]^{-1}\in L(X,X_0)$. As a consequence, $R(\lambda,G'(0))=[\lambda I-G'(0)]^{-1}$ is a compact
  linear operator in $X$, so that its spectrum contains only eigenvalues. It follows that $\sigma(G'(0))$ also contains only eigenvalues. Next, for
  every $s\geqslant 0$ let $H^s({\mathbf{S}}^2)$ be the standard Sobolev space on ${\mathbf{S}}^2$ with index $s$. We know that $H^s({\mathbf{S}}^2)$ has
  an equivalent norm $\|\xi\|_{H^s({\mathbf{S}}^2)}=\displaystyle\Big[\sum_{k=0}^{\infty}\sum_{l=1}^{2k+1}(1\!+\!\lambda_k)^s|c_{kl}|^2\Big]^{1/2}$,
  if $\xi\in H^s({\mathbf{S}}^2)$ has the expression (3.8) as an element of $L^2({\mathbf{S}}^2)$. Since $C^{\infty}({\mathbf{S}}^2)=\displaystyle
  \bigcap_{s\geqslant 0}H^s({\mathbf{S}}^2)$ and it is dense in every $H^s({\mathbf{S}}^2)$, $s\geqslant 0$, Lemma 3.4 shows that the operator $G'(0):
  \dot{C}^{m+3+\mu}({\mathbf{S}}^2)\to\dot{C}^{m+\mu}({\mathbf{S}}^2)$ can be uniquely extended into a bounded linear operator from $H^3({\mathbf{S}}^2)$
  to $L^2({\mathbf{S}}^2)$, and after extension we have the relation $\sigma(G'(0))=\{a_k(\gamma):k=0,1,2,\cdots\}$. Moreover, the eigenspace
  corresponding to the eigenvalue $a_k(\gamma)$ is ${\rm span}\{Y_{kl}(\omega):l=1,2,\cdots,2k\!+\!1\}\subseteq C^{\infty}({\mathbf{S}}^2)$, i.e, all
  eigenfunctions are smooth. Hence, since if $\xi\in\dot{C}^{m+3+\mu}({\mathbf{S}}^2)$ is an eigenvector of the operator $G'(0):
  \dot{C}^{m+3+\mu}({\mathbf{S}}^2)\to\dot{C}^{m+\mu}({\mathbf{S}}^2)$ then it is also an eigenvector of the operator $G'(0):H^3({\mathbf{S}}^2)\to
  L^2({\mathbf{S}}^2)$, we obtain the desired assertion. $\quad\Box$
\medskip

  It is clear that $a_0(\gamma)$ and $a_1(\gamma)$ are independent of $\gamma$. Hence we re-denote them as $a_0$ and $a_1$, respectively, i.e.,
$$
  a_0=g(1)R_s-aR_s\sigma_s'(R_s)\bar{v}_0'(R_s), \quad  a_1=g(1)R_s-aR_s\sigma_s'(R_s)\bar{v}_1'(R_s).
$$
   For $k\geqslant 2$ we denote
\begin{equation}
  \gamma_k=\frac{2R_s^3}{k(k-1)(k+2)}[g(1)-a\sigma_s'(R_s)\bar{v}_k'(R_s)].
\end{equation}
  Then from (3.13) we have
\begin{equation}
  a_k(\gamma)=-\frac{1}{2R_s^2}k(k-1)(k+2)(\gamma-\gamma_k), \qquad k=2,3,\cdots.
\end{equation}
  We shall prove $a_0<0$, $a_1=0$ and $\gamma_k>0$ for $k\geqslant 2$. For this purpose we need the following lemma:
\medskip

  {\bf Lemma 3.6}\ \ {\em For the solution of the problem $(3.9)$ we have the following assertions:

  $(1)$\ \ $0<\bar{u}_k(r)<1$ and $\bar{u}_k'(r)>0$ for $ K_s<r<R_s$.

  $(2)$\ \ If $k>l$ then $\bar{u}_k(r)>\bar{u}_l(r)$ for $ K_s<r<R_s$, and
  $\bar{u}_k'( K_s)\geqslant\bar{u}_l'( K_s)$, $\bar{u}_k'(R_s)\leqslant\bar{u}_l'(R_s)$.

  $(3)$\ \ If $k>l$ then $\bar{u}_k(r)(r/R_s)^{k+1}<\bar{u}_l(r)(r/R_s)^{k+1}$ for $ K_s<r<R_s$.}
\medskip

  {\em Proof.}\ \ The assertion $0<\bar{u}_k(r)<1$ for $ K_s<r<R_s$ is an immediate consequence of the maximum principle. Note that this
  assertion joint with the boundary value conditions $\bar{u}_k( K_s)=0$ and $\bar{u}_k(R_s)=1$ implies that $\bar{u}_k'( K_s)
  \geqslant0$ and $\bar{u}_k'(R_s)\geqslant 0$. Next we let $w_k(r)=\bar{u}_k'(r)$. A simple computation shows that $w_k$ satisfies the following
  equation:
$$
  w_k''(r)+\frac{2(k\!+\!1)}{r}w_k'(r)-\Big(\frac{2(k\!+\!1)}{r^2}+1\Big)w_k(r)=0 \quad  \mbox{for}\;\;  K_s<r<R_s.
$$
  Since $w_k( K_s)\geqslant0$ and $w_k(R_s)\geqslant 0$, again by the maximum principle we see that $w_k(r)>0$ for $ K_s<r<R_s$. This
  proves the assertion (1). From the property $\bar{u}_k'(r)>0$ for $ K_s<r<R_s$ it follows that if $k>l$ then
$$
  \bar{u}_k''(r)+\frac{2(l\!+\!1)}{r}\bar{u}_k'(r)-\bar{u}_k(r)<0 \quad  \mbox{for}\;\;  K_s<r<R_s.
$$
  Hence by the maximum principle we obtain $\bar{u}_k(r)>\bar{u}_l(r)$ for $ K_s<r<R_s$ and $k>l$, which easily implies that $\bar{u}_k'( K_s)
  \geqslant\bar{u}_l'( K_s)$ and $\bar{u}_k'(R_s)\leqslant\bar{u}_l'(R_s)$ for $k>l$. This proves the assertion (2). Finally we
  let $z_k(r)=\bar{u}_k(r)(r/R_s)^{k+1}$. It can be easily seen that $z_k(r)$ is a solution of the following problem:
$$
  \left\{
\begin{array}{l}
  \displaystyle z_k''(r)-\Big(\frac{k(k\!+\!1)}{r^2}+1\Big)z_k(r)=0 \quad  \mbox{for}\;\;  K_s<r<R_s,\\ [0.2cm]
   z_k( K_s)=0,   \quad  z_k(R_s)=1.
\end{array}
\right.
$$
  Since $z_k>0$ for $ K_s<r<R_s$, a similar argument as in the proof of the assertion (2) shows that if $k>l$ then $z_k(r)<z_l(r)$ for  for
  $ K_s<r<R_s$. This proves the assertion (3) and completes the proof of the lemma. $\quad\Box$
\medskip

  {\bf Lemma 3.7}\ \ {\em We have the following assertions:

  $(1)$\ \ $a_0<0$ and $a_1=0$.

  $(2)$\ \ $\gamma_k>0$ for all $k=2,3,\cdots$, and $\gamma_k\sim 2R_s^3g(1)k^{-3}$ as $k\to\infty$.}
\medskip

  {\em Proof.}\ \ From (3.9) we have
$$
  \bar{u}_k'(R_s)=\bar{u}_k'( K_s)\Big(\frac{ K_s}{R_s}\Big)^{2(k+1)}+
  \int_{ K_s}^{R_s}\bar{u}_k(\tau)\Big(\frac{\tau}{R_s}\Big)^{2(k+1)}d\tau.
$$
  Hence the relation (3.12) can be rewritten as follows:
$$
  \bar{v}_k'(R_s)=\frac{\hat{\sigma}-\tilde{\sigma}}{\hat{\sigma}}\bar{u}_k'(R_s)
  +\frac{\tilde{\sigma}}{\hat{\sigma}}\int_{ K_s}^{R_s}\bar{u}_k(\tau)\Big(\frac{\tau}{R_s}\Big)^{2(k+1)}d\tau.
$$
  Using this relation and the assertions (2) and (3) of Lemma 3.6, we conclude that
$$
  \bar{v}_k'(R_s)<\bar{v}_l'(R_s) \quad \mbox{if}\;\; k>l.
$$
  Hence, all the desired assertions will follow if we prove that $a_1=0$ and $\bar{v}_k'(R_s)=O(1/k)$ as $k\to\infty$. The proof that
  $\bar{v}_k'(R_s)=O(1/k)$ as $k\to\infty$ is easy and is omitted. In what follows we prove $a_1=0$.

  It is easy to check
\begin{equation}
  \bar{u}_1(r)=\frac{R_s\sigma_s'(r)}{r\sigma_s'(R_s)} \quad \mbox{for}\;\,  K_s\leqslant r\leqslant R_s.
\end{equation}
  In what follows we prove
\begin{equation}
  \bar{v}_1(r)=-\frac{R_s\pi_s'(r)}{ar\sigma_s'(R_s)} \quad \mbox{for}\;\,  K_s\leqslant r\leqslant R_s.
\end{equation}
  We use the notation $q(r)$ to denote the function on the right-hand side of the above relation. To prove the above relation, we only need to
  show that $q(r)$ is a solution of the following problem:
$$
  \left\{
\begin{array}{l}
  \displaystyle q''(r)+\frac{4}{r}q'(r)=\bar{u}_1(r) \quad  \mbox{for}\;\;  K_s<r<R_s,\\ [0.2cm]
  \displaystyle  q(R_s)=0,   \quad  q'(K_s^+)=\frac{\hat{\sigma}-\tilde{\sigma}}{\hat{\sigma}}\bar{u}_1'(K_s^+).
\end{array}
\right.
$$
  The equation in the first line is easy to check, and the boundary value condition $q(R_s)=0$ is clear. Since
$$
  q'(r)=-\frac{R_s}{a\sigma_s'(R_s)}\Big[\frac{\pi_s''(r)}{r}-\frac{\pi_s'(r)}{r^2}\Big]
  =\frac{R_s}{a\sigma_s'(R_s)}\Big[\frac{3\pi_s'(r)}{r^2}+\frac{g(\sigma_s(r))}{r}\Big],
$$
  we have
$$
  q'( K_s^+)=\frac{R_s}{a\sigma_s'(R_s)}\Big[\frac{3\pi_s'( K_s)}{ K_s^2}+\frac{g(\hat{\sigma}^+)}{ K_s}\Big]
  =\frac{3R_s\pi_s'(K_s)}{a\sigma_s'(R_s)K_s^2}.
$$
  From the equation $\Delta\pi_s=-b$ (for $r<K_s$) we see that $\pi_s'(K_s)=(1/3)bK_s=(a/3)(\hat{\sigma}-\tilde{\sigma})K_s$, and it is clear
  that $\bar{u}_1'(K_s^+)=R_s\sigma_s''(K_s^+)/K_s\sigma_s'(R_s)=\hat{\sigma}R_s/K_s\sigma_s'(R_s)$. Combining these relations, we see that the
  the boundary value condition $q'(K_s^+)=\frac{\hat{\sigma}-\tilde{\sigma}}{\hat{\sigma}}\bar{u}_1'( K_s^+)$ is also satisfied. Hence (3.19) is
  true. The assertion $a_1=0$ is an immediate consequence of the relation (3.19) and the fact that $\pi_s''(R_s)=-g(1)$. $\quad\Box$
\medskip

  {\bf Lemma 3.8}\ \ {\em Let $\gamma^*=\max\{\gamma_k:k=2,3,\cdots\}$. The following assertions hold:

  $(1)$\ \ If $\gamma>\gamma^*$ then $\sup\{{\rm Re}\lambda:\lambda\in\sigma(G'(0))\backslash\{0\}\}<0$ and ${\rm Ker}\,G'(0)=
  {\rm span}\{Y_{11}(\omega),Y_{12}(\omega),Y_{13}(\omega)\}$, so that $\dim{\rm Ker}\,G'(0)=3$. If instead $0<\gamma<\gamma^*$ then
  $\sup\{{\rm Re}\lambda:\lambda\in\sigma(G'(0))\}>0$.

  $(2)$\ \ Let $\gamma>\gamma^*$. Then ${\rm Range}\,G'(0)$ is closed, and $X={\rm Ker}\,G'(0)\oplus{\rm Range}\,G'(0)$.}
\medskip

  {\em Proof.}\ \ Assertions in (1) are immediate consequences of Corollary 3.5, the expression (3.15) and Lemma 3.7.
  To prove the assertion (2), we note that $G'(0)=\gamma F'(0)+F_0'(0)$. From \cite{EscSim2} we know that $F'(0)$ is
  a third-order elliptic pseudo-differential operator on the sphere ${\mathbf{S}}^2$, and Lemma 2.5 shows that $F_0'(0)$
  is a lower-order perturbation. It follows that standard $C^{\mu}$-estimates work for $G'(0)$ and, consequently, the
  Fredholm alteration principle applies to it, by a similar argument as in the proof of Lemma 2.5. Hence the assertion
  (2) follows. $\quad\Box$
\medskip

  {\em Proof of Theorem 1.2}:\ \ From (3.15) and Corollary 3.5 (as well as the fact that $a_0<0$) we see that if
  $\gamma>\gamma^*$ then $\sup\{\mbox{Re}\lambda:\lambda\in\sigma(G'(0))\}<0$. This fact and Lemmas 3.1, 3.2, 3.3, 3.7, 3.8
  show that Theorem 1.1 of \cite{Cui5} (with $N=1$) applies to the equation $(3.3)_1$. Hence, by applying Theorem 1.1 of
  \cite{Cui5} we see that if $\gamma>\gamma^*$ then the assertions (1)--(3) of Theorem 1.2 hold. If on the other hand
  $\gamma<\gamma^*$ then from (3.15) we see that there exists integer $k\geqslant 2$ such that $a_k(\gamma)>0$, so that
  $\sup\{\mbox{Re}\lambda:\lambda\in\sigma(G'(0))\}>0$. It follows from linearized instability criterion for parabolic
  equations in Banach spaces (i.e. Theorem 9.1.3 of \cite{Lun2}) we obtain the last assertion of Theorem 1.2. This
  proves Theorem 1.2. $\quad\Box$
\medskip


\end{document}